\documentclass[hidelinks,onefignum,onetabnum]{siamart251216}

\usepackage{lipsum}
\usepackage{amsfonts}
\usepackage{graphicx}
\usepackage{epstopdf}
\usepackage{amsmath,amssymb}
\usepackage{graphicx}
\usepackage{subcaption}
\usepackage{algorithm,todonotes}
\usepackage{algpseudocode}
\usepackage{mathtools}
\usepackage{mathrsfs}
\usepackage{color,xcolor}
\newcommand{\jn} [1]{{\color{black} #1}}
\newcommand{\dn}[1]{{\color{black}#1}}
\usepackage{hyperref}
\hypersetup{colorlinks=true, urlcolor=blue, citecolor=blue, linkcolor=blue}
\ifpdf
\DeclareGraphicsExtensions{.eps,.pdf,.png,.jpg}
\else
  \DeclareGraphicsExtensions{.eps}
\fi

\newcommand{\Loss}{\mathcal{L}}
\newsiamremark{remark}{Remark}
\newsiamremark{hypothesis}{Hypothesis}
\crefname{hypothesis}{Hypothesis}{Hypotheses}
\newsiamthm{claim}{Claim}
\newsiamremark{fact}{Fact}
\crefname{fact}{Fact}{Facts}

\headers{New applications for optimal adaptivity}{M. Feischl, J. Nick, D. Niederkofler}

\def\R{{\mathbb R}}

\def\N{{\mathbb N}}

\def\SS{{\mathcal S}}

\def\TT{{\mathcal T}}

\def\tend{{t_{\rm end}}}

\def\bf{{\boldsymbol{f}}}

\def\norm#1#2{\|#1\|_{#2}}

\def\eps{\varepsilon}

\def\normL2#1#2{\|#1\|_{L^2(#2)}}

\def\revision#1{{\color{black}{#1}}}

\def\T{\mathbb T}

\title{New applications for optimal adaptivity\thanks{Submitted to the editors DATE.
}}

\author{Michael Feischl\thanks{TU Wien,
  \email{michael.feischl@asc.tuwien.ac.at}} \and Jörg Nick
  \thanks{University of Birmingham,
  \email{j.nick@bham.ac.uk}} 
\and David Niederkofler\thanks{TU Wien,
  \email{david.niederkofler@asc.tuwien.ac.at}}}

\usepackage{amsopn}

\begin{document}

\maketitle

\begin{abstract}
    \dn{We introduce \revision{an explicit time-stepping scheme} for ODEs, together with a residual-based a posteriori error estimator.}
    \dn{Using this estimator, }we show \revision{a new approach to proving optimality and linear reduction of adaptive algorithms for time-dependent problems} that works in $H^1(0,\tend)$ and $W^{1,\infty}(0,\tend)$ norms, \revision{does not require \emph{quasi-orthogonality}, and poses} no restrictions on the marking parameter $\theta$. The key idea is to use the Lipschitz continuity of the semi-discretized operator to obtain \revision{a stronger form of estimator reduction}. \dn{Numerical experiments illustrate the benefits of this explicit scheme coupled with adaptive refinement over uniform time stepping for classical evolution problems, regularized dynamical approximation, and neural ODEs.}
\end{abstract}

\section{Introduction}
This work establishes a class of model problems and corresponding adaptive methods that showcase \dn{the following contributions}:
\begin{itemize}
    \item Optimal adaptive convergence rates without the use of \emph{quasi-orthogonality} as a theoretical tool.
\end{itemize}
Even the first optimality proofs in~\cite{stevenson07, ckns, bdd} use orthogonality in the form of \revision{the Pythagorean identity}. For increasingly involved model problems, several extensions of this orthogonality have been established under the name of quasi-orthogonality, e.g.,~\cite{nonconf, ccgs13,CCnonlin} for non-conforming elements and~\cite{nonsymm,feischl2022} for non-symmetric and indefinite problems. However, all versions of this quasi-orthogonality are technically involved and pose severe restrictions on the problem class, particularly in terms of nonlinearities. The use of explicit time stepping together with the Lipschitz continuity of the problem allows us to remove any form of quasi-orthogonality from the proof and thus also to remove the restrictions on the \revision{nonlinearity} for $p=1$. 
\begin{itemize}
    \item Optimal convergence rates in $W^{1,\infty}$ with maximum marking.
\end{itemize}
Except for convergence results in the $L^2$-norm for an $H^1$-problem in~\cite{stevensonl2}, all optimality proofs work in the natural norm of the problem.
\dn{Moreover, optimality proofs for adaptive algorithms have predominantly
relied on D\"orfler marking~\cite{d1996}, with exceptions including
instance optimality for a modified maximum strategy~\cite{instanceopt}
and the recent proof of rate optimality for standard maximum marking
for the Poisson problem~\cite{li2026rateoptimalityafemmaximum}.}
The simplified optimality analysis that comes with the removal of quasi-orthogonality has the added benefit that it also works in the $W^{1,\infty}$-norm. It turns out that maximum marking is the natural candidate for this norm and hence many of the well-known arguments from the energy norm case transfer almost verbatim.
\begin{itemize}
    \item Provable optimality of the adaptive algorithm without restrictions on the marking parameter $0<\theta<1$.
\end{itemize}
The classical proof strategy for optimal convergence \revision{requires choosing the marking parameter} $\theta$ \revision{in D\"orfler marking sufficiently small}~\cite{axioms}. Essentially, this means that the adaptive refinement must force aggressive local refinement to be optimal. This seems to be largely a theoretical issue of the proof as almost no practical examples are known where this behavior can be seen. The work~\cite{thetaopt} even shows that the theoretical framework in~\cite{axioms} is not strong enough to prove optimality \revision{without restrictions} \dn{\revision{on the marking parameter.} Recently, the work \cite{diening2026optimalafemconvergencebulk} was able to lift the restriction on the marking parameter for the Poisson problem.}
The key new tool in the present work is a strengthened version of estimator reduction (see~\cite{axioms} for an overview on this proof technique), which also removes the smallness restrictions on the marking parameter, both in the D\"orfler marking and in the maximum marking. \dn{\revision{The tradeoff is the need for a fine initial mesh.}}
\begin{itemize}
    \item \dn{Application to \revision{neural ODEs and dynamical parametric approximation}.}
\end{itemize}
Neural ODEs replace the classical layer structure of feed-forward neural networks \revision{with an ODE that has to be solved} for training and inference of the network. \revision{The authors of the seminal paper}~\cite{neuralODE} \revision{already noticed that adaptive time stepping is necessary} to solve the sometimes unstable equations. Our explicit adaptive schemes are versatile enough to be applied in this setting and provide a clear advantage over standard uniform training of neural ODEs.
\dn{A further application is regularized dynamical parametric
approximation~\cite{regdyn,LN25}, where the solution of an evolution
equation is approximated by a nonlinear parametrization
$u(t)=\Phi(q(t))$ with time-dependent parameters $q(t)$. For example,
$\Phi$ may represent a neural network whose weights and biases evolve
in time. We extend our adaptive
framework to this setting and \revision{provide numerical experiments} with a neural-network
parametrization of a transport problem to illustrate the benefits of
adaptive time stepping.} 
\section{An explicit discretization method}\noindent
\revision{We first introduce the discretization and residual estimator that form the basis of the adaptive analysis.}
We consider autonomous initial value problems of the form
\begin{align}
\label{eq:ode}
\begin{split}
    \partial_t y(t) &= F(y(t)), \qquad \text{for } t \in [t_0, \tend]\\
    y(t_0) &= y_0,
\end{split}
\end{align}
for a Lipschitz continuous right-hand side $F:\mathbb{R}^d \to \mathbb{R}^d$. The time interval $[t_0, \tend]$ is discretized into a finite number of subintervals $\TT=\{[t_0,t_1],\ldots, [t_{n-1},\tend]\}$. We use the standard spline spaces  $P^p(\TT)=\{v: v_{|T} \in P^p(T)\}$ and $\SS^p(\TT)=P^p(\TT) \cap C([t_0,\tend])$. The discretization scheme we propose is similar to Petrov-Galerkin schemes in time, but with one significant difference: For $T=[t_i,t_{i+1}]$, we introduce the local operator $I^F_T: \mathbb{R}^d \to P^{p-1}(T)$, mapping some initial value on the interval $T$ to a \revision{lower-order} polynomial approximation of $y$ on the same interval. \jn{As a result, we obtain a \revision{lower-order} explicit time-stepping scheme.} Using this \revision{lower-order} scheme, we aim to build a \revision{higher-order} approximation $y_\TT \in \SS^p(\TT)$ satisfying for all $T \in \TT$ and $v \in P^{p-1}(\TT)$,
\begin{align}
\label{eq:ode_disc}
\begin{split}
    y_\TT(t_0)=y_0\quad\text{and}\quad
    \int_T \Big( \partial_t y_\TT - F\big(I^F_T(y_\TT(t_i))\big) \Big) \cdot v \,\mathrm{d}t =0.
\end{split}
\end{align}
It is apparent that this yields an explicit time-stepping scheme. For error estimation, we use estimators based on the residual $R=\partial_t y_\TT-F(y_\TT)$, as a Gronwall argument gives the following estimate.
\begin{lemma}\label{lem:rel}
    Let $y$ solve \eqref{eq:ode} and $y_\TT \in \SS^p(\TT)$, with $y_\TT(t_0)=y_0$. For the residual $R=\partial_t y_\TT-F(y_\TT)$ it holds that
    \begin{align*}
        \norm{y-y_\TT}{H^1(t_0,\tend)} &\leq  C_{\rm rel} \norm{R}{L^2(t_0,\tend)},\\
        \norm{y-y_\TT}{W^{1,\infty}(t_0,\tend)} &\leq  C_{\rm rel} \norm{R}{L^\infty(t_0,\tend)},
    \end{align*}
    with a constant $C_{\rm rel}>1$.
\end{lemma}
Concretely, the estimator will be of the form
     $\eta_{\TT}^2= \sum_{T \in \TT} \eta_\TT(T)^2.$
With this, the standard adaptive loop in Algorithm~\ref{alg:adaptive} can be formulated. {\color{black} Related residual-based error estimators for explicit Runge--Kutta methods are studied by Ranocha and Giesselmann~\cite{explicit_a_posteriori}, who use the residual of a suitable reconstruction to select time-step sizes sequentially during forward integration. In contrast, our adaptive loop uses the local error indicators from the entire time interval to mark subintervals for refinement and then recomputes the discrete solution on the refined mesh.}
\begin{algorithm}[h]
    \caption{Adaptive Algorithm}
    \label{alg:adaptive}
    \begin{algorithmic}[1]
       \item \textbf{Input:} Initial mesh $\TT_0$, parameter $0<\theta<1$.
        \For{ $\ell=0,1,2,\ldots$}:
        \State Compute $y_{\TT_\ell}$ from \eqref{eq:ode_disc}.
        \State Compute error estimates $\eta_{\TT_\ell}(T)$ for all $T \in \TT_\ell$.
        \State Find a set $\mathcal{M}_\ell$ of minimal cardinality such that
        \begin{align*}
            \sum_{T \in \mathcal{M}_\ell} \eta_{\TT_\ell}(T)^2 \geq \theta \sum_{T \in \TT_\ell} \eta_{\TT_\ell}(T)^2.
        \end{align*}
        \State Bisect all $T \in \mathcal{M}_\ell$ to obtain a new mesh $\TT_{\ell+1}$.
        \EndFor
    \end{algorithmic}
\end{algorithm}
The main result of this work is to prove \revision{so-called} rate-optimality for the adaptive algorithm, in the sense of~\cite{axioms}. To that end, we introduce the set $\T$ of all  meshes that can be reached by \revision{repeated bisections}, starting from some initial mesh $\TT_0$.
\begin{theorem}\label{thm:optim}
\revision{Assume the setting of Section~}\ref{sec:low_order} or~\ref{sec:high_order}.
    Let $0<\theta<1$ and suppose the initial discretization $\TT_0$ is sufficiently fine. Let $s>0$ be  such that
    \begin{align*}
        C_{\rm best}:=\sup_{N \in \mathbb{N}} \inf_{\substack{\TT \in \T \\ \#\TT-\#\TT_0\leq N}}(N+1)^s \eta_\TT<\infty.
    \end{align*}
    Then, the sequence of meshes $\TT_\ell$ obtained from Algorithm~\ref{alg:adaptive} satisfies 
\begin{align}\label{eq:rates}
        \sup_{\ell \in \mathbb{N}} \eta_{\TT_\ell}(\#\TT_\ell-\#\TT_0+1)^s \leq C_{\rm opt} C_{\rm best}.
    \end{align}
\end{theorem}
\subsection{\revision{Lowest-order} discretization $p=1$}\noindent \label{sec:low_order}
To obtain an approximation scheme for $p=1$, we use the mapping $I_T^F(z)=z\in P^0(T)$. Evaluating \eqref{eq:ode_disc} yields the \revision{well-known} explicit \revision{time-stepping method}
\begin{align}\label{eq:exp_euler}
    y_\TT(t)=y_\TT(t_i)+(t-t_i)F(y_\TT(t_i)), \quad t \in [t_i,t_{i+1}].
\end{align}
In this setting, the residual is given by
\begin{align*}
    R(t)=F(y_\TT(t_i))-F(y_\TT(t)), \quad t \in [t_i,t_{i+1}].
\end{align*}
The Lipschitz condition on $F$ allows us to derive \revision{the estimator that will be used} in this setting. For $T=[t_i,t_{i+1}]$ one gets
\begin{align}\label{eq:estimator}
\eta_\TT(T) \coloneq L |T| \norm{\partial_t y_\TT}{L^2(T)}=  L|T|\norm{F(y_\TT(t_i))}{L^2(T)} \geq
    \norm{R}{L^2(T)} .
\end{align}
We aim to prove Theorem~\ref{thm:optim} for this pair of discretization and estimator. \revision{In contrast to the classical strategy outlined in \cite{axioms},} we demonstrate \textit{stronger} stability and reduction \revision{estimates that are sufficient} to prove R-linear convergence of the estimator.
\begin{lemma}\label{lem:stab_1}
    The error estimator defined in \eqref{eq:estimator} satisfies \textit{strong} stability in the sense that for a refinement $\widehat\TT$ of $\TT$ it holds
    \begin{align*}
        \Bigg| \Big(\sum_{T \in \TT \cap \widehat\TT}\eta_{\TT}(T)^2 \Big)^{1/2} &-\Big(\sum_{T \in \TT \cap \widehat\TT} \eta_{\widehat\TT}(T)^2\Big)^{1/2} \Bigg| \\
 &\leq  C_{\rm stab} \norm{h_\TT}{L^\infty(t_0,\tend)}\norm{y_{\TT}-y_{\widehat\TT}}{H^1(t_0,\tend)},
    \end{align*}
     where $y_\TT$ and $y_{\widehat{\TT}}$ are discrete solutions on the respective meshes, $h_\TT(t)= |T|$ for $t \in \text{int}(T)$ is the mesh-size function and $C_{\rm stab}>1$.
\end{lemma}
\begin{proof}
    Triangle inequalities show
    \begin{align*}
        \Bigg| \Big(\sum_{T \in \TT \cap \widehat\TT}\eta_{\TT}(T)^2 \Big)^{1/2} &-\Big(\sum_{T \in \TT \cap \widehat\TT} \eta_{\widehat\TT}(T)^2\Big)^{1/2} \Bigg|^2
        \\
 &\leq  \sum_{T \in \TT \cap \widehat\TT} L^2 |T|^2 \norm{\partial_t y_\TT-\partial_t y_{\widehat\TT}}{L^2(T)}^2\\
 &\leq  L^2 \norm{h_\TT}{L^\infty(t_0,\tend)}^2 \norm{\partial_ty_\TT-\partial_t y_{\widehat\TT}}{L^2(t_0,\tend)}^2.
    \end{align*}
\end{proof}
\revision{A similar estimate on refined elements yields} \textit{strong} reduction.
\begin{lemma}\label{lem:red_1}
    The error estimator defined in \eqref{eq:estimator} satisfies \textit{strong} reduction in the sense that for a refinement $\widehat\TT$ of $\TT$ it holds
    \begin{align*}
    \sum_{T \in \widehat\TT \setminus \TT} \eta_{\widehat\TT}(T)^2 &\leq  q \sum_{T \in  \TT\setminus \widehat\TT} \eta_{\TT}(T)^2 \\
 &\quad + C_{\rm red}\norm{h_\TT}{L^\infty(t_0,\tend)}^2\norm{y_{\TT}-y_{\widehat\TT}}{H^1(t_0,\tend)}^2,
\end{align*}
    for $0<q<1$ and $C_{\rm red}>1$.
\end{lemma}
\begin{proof}
    Let $T \in \TT \setminus \widehat\TT$ and  $\{\widehat{T}_i\}_{i=1}^n \subset \widehat\TT \setminus \TT$ such that $T=\bigcup_{i=1}^n \widehat{T}_i $. It holds
    \begin{align*}
       \sum_{i=1}^n \eta_{\widehat\TT}(\widehat{T}_i)^2 &=  \sum_{i=1}^n L^2|\widehat{T}_i|^2\norm{\partial_t y_{\widehat\TT}}{L^2(\widehat{T}_i)}^2\\
 &\leq   L^2\frac{|T|^2}{2}\norm{\partial_t y_{\TT}}{L^2(T)}^2 + L^2\frac{|T|^2}{2}\norm{\partial_t y_{\TT}-\partial_t y_{\widehat\TT}}{L^2(T)}^2 \\
 &\leq  \frac12 \eta_\TT(T)^2 + L^2\frac{|T|^2}{2}\norm{\partial_t y_{\TT}-\partial_t y_{\widehat\TT}}{L^2(T)}^2,
    \end{align*}
    \revision{since} $|\widehat{T}_i|\leq \frac12 |T|$. \revision{Summing over all $T \in \TT \setminus \widehat\TT$} yields the result.
\end{proof}
 Following the standard theory for adaptive mesh refinement, one would continue to prove the \revision{so-called} quasi-orthogonality \revision{estimate in order to conclude} R-linear convergence of the estimator as outlined in \cite{axioms,feischl2022}. \revision{In this work, we use} the stronger stability and reduction estimates to conclude linear convergence directly, under an additional assumption on the initial mesh.
\begin{lemma}
    \label{lem:rlinearconv}
    Let $\TT_0, \TT_{1},\ldots$ be meshes generated by Algorithm~\ref{alg:adaptive} with \revision{corresponding discrete solutions} $y_{0}, y_{1},\ldots$ and estimators $\eta_{0}, \eta_{1},\ldots$. Suppose the initial discretization $\TT_0$ is fine enough. Furthermore, suppose that the estimator satisfies strong stability and reduction as outlined in Lemmas~\ref{lem:stab_1} and~\ref{lem:red_1}. Then there is a $\kappa\in(0,1)$ 
    \revision{such that, for all $k,\ell \in \N$,}
    \begin{align*}
        \eta^2_{\ell+k} \leq \kappa^k \eta^2_\ell.
    \end{align*}
\end{lemma}
\begin{proof}
    Let $\ell \in \N$. \revision{Using the stability estimate in Lemma~\ref{lem:stab_1}, the reduction estimate in Lemma~\ref{lem:red_1}, and Young's inequality, we obtain}
        \begin{align*}
            \eta_{\ell+1}^2 &=  \sum_{T \in \TT_{\ell+1} \setminus \TT_\ell} \eta_{\ell+1}(T)^2 + \sum_{T \in \TT_{\ell+1} \cap \TT_{\ell}} \eta_{\ell+1}(T)^2\\
 &\leq  q \sum_{T \in \TT_{\ell} \setminus \TT_{\ell+1}} \eta_{\ell}(T)^2+ (1+\delta) \sum_{T \in \TT_{\ell+1} \cap \TT_{\ell}} \eta_{\ell}(T)^2\\
 &+ C_{\rm est} \norm{h_{\TT_\ell}}{L^\infty(t_0,\tend)}^2\norm{y_{\ell}-y_{\ell+1}}{H^1(t_0,\tend)}^2
        \end{align*}
        with $C_{\rm est}=C_{\rm red}+(1+\delta^{-1})C^2_{\rm stab}$, and $\delta>0$ arbitrary. Using the marking strategy with $\theta \in (0,1)$ yields
        \begin{align}
        \label{eq:est_reduction}
             \eta_{\ell+1}^2 \leq ((1+\delta)(1-\theta)+\theta q) \eta^2_\ell + C_{\rm est} \norm{h_{\TT_\ell}}{L^\infty(t_0,\tend)}^2\norm{y_{\ell}-y_{\ell+1}}{H^1(t_0,\tend)}^2.
        \end{align}
        \revision{Choosing $\delta$ small enough ensures that} $q_{\rm red}:=((1+\delta)(1-\theta)+\theta q)<1$. \revision{The estimate \eqref{eq:est_reduction} is called estimator reduction} in \cite{axioms}. Since we have the reliability estimate
        \begin{align*}
            \norm{y_{\ell}-y_{\ell+1}}{H^1(t_0,\tend)}^2 \leq 2 C_{\rm rel}^2 \Big(  \eta^2_\ell+\eta^2_{\ell+1} \Big),
        \end{align*}
        we get for $C=2C_{\rm rel}^2 C_{\rm est}$ that
        \begin{align*}
            \eta_{\ell+1}^2 (1-C\norm{h_{\TT_\ell}}{L^\infty(t_0,\tend)}^2) \leq \eta_{\ell}^2 (q_{\rm red}+C\norm{h_{\TT_\ell}}{L^\infty(t_0,\tend)}^2),
        \end{align*}
        which yields under the assumption of a fine initial mesh such that $C \norm{h_{\TT_\ell}}{L^\infty(t_0,\tend)}^2\leq C\norm{h_{\TT_0}}{L^\infty(t_0,\tend)}^2<1$ that
        \begin{align*}
            \eta_{\ell+1}^2 \leq \frac{q_{\rm red}+C\norm{h_{\TT_0}}{L^\infty(t_0,\tend)}^2}{1-C\norm{h_{\TT_0}}{L^\infty(t_0,\tend)}^2} \eta_{\ell}^2.
        \end{align*}
        \revision{Choosing an even finer initial discretization $\TT_0$ gives} 
        \begin{align*}
            \kappa \coloneq \frac{q_{\rm red}+C\norm{h_{\TT_0}}{L^\infty(t_0,\tend)}^2}{1-C\norm{h_{\TT_0}}{L^\infty(t_0,\tend)}^2}<1,
        \end{align*} which concludes the proof.
\end{proof}
In order to prove optimal convergence rates of the algorithm, one still needs another crucial property of the error estimator, namely a discrete version of the reliability estimate. To obtain such a result, we state an auxiliary lemma.
\begin{lemma}\label{lem:help_disc_rel}
    Let $\widehat\TT$ be a refinement of $\TT$. Suppose that for all unrefined elements $T \in \TT \cap \widehat\TT$, the approximation scheme \revision{satisfies the following estimate:}
    \begin{align}\label{eq:initial_esti}
        |y_{\widehat\TT}(t)-y_\TT(t)| \leq (1+C |T|)|y_{\widehat\TT}(t_i)-y_\TT(t_i)|, \quad t \in T=[t_i,t_{i+1}],
    \end{align}
    where $C>1$ is some constant independent of $\TT$ and $\widehat \TT$.
    Furthermore, assume that the initial mesh $\TT_0$ is fine enough. \revision{Then, for all $t \in [t_0,\tend]$,}
    \begin{align}
        |y_{\widehat\TT}(t)-y_\TT(t)| \lesssim \sum_{T \in \widehat\TT \setminus \TT} \int_T \Big(|R(s)|+|\widehat{R}(s)|\Big)\,\mathrm{d}s,
    \end{align}
    where $R=\partial_t y_\TT-F(y_\TT)$ and $\widehat{R}=\partial_t y_{\widehat\TT}-F(y_{\widehat\TT})$ denote the corresponding residuals.
\end{lemma}
\begin{proof}
    Let $[t_i, t_{i+1}]=T \in \widehat\TT \setminus \TT$. With the residuals $\partial_t y_\TT-F(y_\TT):= R$ and $\partial_t y_{\widehat\TT}-F(y_{\widehat\TT}):= \widehat R$ it holds
\begin{align*}
    \partial_t y_{\widehat\TT} -\partial_t y_{\TT}= F(y_{\widehat\TT})-F(y_\TT)+\widehat R(t)-R(t).
\end{align*}
Integrating from $t_i$ to $t \in T$ yields
\begin{align*}
 &|y_{\widehat\TT}(t)-y_{\TT}(t)| \\
 &\leq  |y_{\widehat\TT}(t_i)-y_{\TT}(t_i)|+ \int_T |F(y_{\widehat\TT}(s))-F(y_{\TT}(s))| \,\mathrm{d}s
    +\int_T \Big(|R(s)|+|\widehat R(s)| \Big) \,\mathrm{d}s,\\
 &\leq  |y_{\widehat\TT}(t_i)-y_{\TT}(t_i)|+ L|T| \norm{y_{\widehat\TT}-y_{\TT}}{L^\infty(T)}
    +\int_T \Big(|R(s)|+|\widehat R(s)| \Big) \,\mathrm{d}s,
\end{align*}
where we used the Lipschitz condition in the second step. This results in
\begin{align}\label{eq:ref_esti}
    \begin{split}
        \norm{y_{\widehat\TT}-y_{\TT}}{L^\infty(T)} &\leq  \frac{1+C|T|}{1-L|T|}\Big(|y_{\widehat\TT}(t_{i})-y_{\TT}(t_{i})| \\
 &\quad + \int_T \Big(|R(s)|+|\widehat R(s)| \Big) \,\mathrm{d}s \Big),
    \end{split}
\end{align}
if $L|T|<1$, which can be enforced by choosing a fine enough initial mesh.
The claim follows by an inductive argument using $y_\TT(t_0)=y_{\widehat\TT}(t_0)$ \revision{together with \eqref{eq:initial_esti} and \eqref{eq:ref_esti} on unrefined and refined elements, respectively.} We number the elements in $\widehat \TT$ such that $T_i:=[t_{i-1},t_i] \in \widehat \TT$ and claim that for $t \in [t_0,t_i]$ and all $i \in \{1,\ldots,\# \widehat \TT\}$ it holds that
\begin{align}\label{eq:ind_hyp}
    \begin{split}
 &|y_{\widehat\TT}(t)-y_\TT(t)| \\
 &\quad\leq \prod_{j=1}^i \frac{1+C|T_j|}{1-L|T_j|} \sum_{j=1 \atop T_j \in \widehat \TT \setminus \TT }^i \int_{T_j} \Big(|R(s)|+|\widehat{R}(s)|\Big) \,\mathrm{d}s.
    \end{split}
\end{align}
\revision{Setting $i=\# \widehat \TT$ yields the result.}
    The base case follows by applying \eqref{eq:initial_esti} and \eqref{eq:ref_esti} for $T=T_1$ using that $y_\TT(t_0)=y_{\widehat\TT}(t_0)$. To prove the induction step for a refined element $T_{i+1} \in \widehat \TT \setminus \TT$ we use \eqref{eq:ref_esti} and insert \eqref{eq:ind_hyp} for $y_{\widehat\TT}(t_{i})-y_{\TT}(t_{i})$ which yields for $t \in T_{i+1}$
    \begin{align*}
        |y_{\widehat\TT}(t)-y_\TT(t)| \leq \frac{1+C|T_{i+1}|}{1-L|T_{i+1}|} &\Bigg( \prod_{j=1}^i \frac{1+C|T_j|}{1-L|T_j|} \sum_{j=1 \atop T_j \in \widehat \TT \setminus \TT }^i \int_{T_j} \Big(|R(s)|+|\widehat{R}(s)|\Big)\\
 &+\int_{T_{i+1}} \Big(|R(s)|+|\widehat R(s)| \Big) \,\mathrm{d}s\Bigg)\\
 &\leq  \prod_{j=1}^{i+1} \frac{1+C|T_j|}{1-L|T_j|} \sum_{j=1 \atop T_j \in \widehat \TT \setminus \TT }^{i+1} \int_{T_j} \Big(|R(s)|+|\widehat{R}(s)|\Big) \,\mathrm{d}s.
    \end{align*}
    \revision{In the case of an unrefined element} $T_{i+1} \in \widehat \TT \cap \TT$ one just uses \eqref{eq:initial_esti} instead of \eqref{eq:ref_esti} and proceeds analogously. This concludes the induction. \dn{Because $L|T_i| \leq L\norm{h_{\TT_0}}{L^\infty(t_0,\tend)}<1$ \revision{the products involved} are bounded independently of $\TT$.} This concludes the proof.
\end{proof}
\revision{Together with the properties of the numerical scheme, Lemma~\ref{lem:help_disc_rel} yields the following discrete reliability estimate.} 
\begin{lemma}\label{lem:discrel}
    Let $\widehat\TT$ be a refinement of $\TT$. Under the assumption of a fine enough initial mesh $\TT_0$, \revision{discrete reliability holds in the sense}
        \begin{align*}
        \norm{y_{\widehat\TT}-y_{\TT}}{H^1(t_0,\tend)}^2 \leq C_{\rm drel}^2 \sum_{T \in \TT \setminus \widehat\TT}\eta_{\TT}(T)^2,
    \end{align*}
    for corresponding discrete solutions $y_{\widehat\TT}, y_\TT$ defined in \eqref{eq:exp_euler} and $C_{\rm drel}>1$, independent of the meshes $\TT$ and $\widehat\TT$.
\end{lemma}
\begin{proof}
 We first consider $[t_i, t_{i+1}]=T \in \TT \cap \widehat\TT$. From the definition of the numerical scheme \eqref{eq:exp_euler} we get that
    \begin{align*}
        |\partial_t y_\TT-\partial_t y_{\widehat\TT}|=\Big|F(y_\TT(t_i))-F(y_{\widehat\TT}(t_i))\Big|,
    \end{align*}
    which yields by using the Lipschitz condition on $F$ \dn{and norm equivalence on \revision{finite-dimensional spaces}, followed by scaling} that
    \begin{align*}
    \norm{\partial_t y_\TT-\partial_t y_{\widehat\TT}}{L^2(T)} \lesssim \norm{y_\TT-y_{\widehat\TT}}{L^2(T)}.
\end{align*}
Furthermore, we also get the estimate \eqref{eq:initial_esti} by using similar arguments. For $[t_i,t_{i+1}]=T \in \TT \setminus \widehat\TT$ we get
\begin{align*}
 &\norm{\partial_t y_\TT-\partial_t y_{\widehat\TT}}{L^2(T)}^2\ \\
 &\lesssim  \norm{\partial_t y_\TT-F(y_\TT)}{L^2(T)}^2+\norm{F(y_\TT)-F(y_{\widehat\TT})}{L^2(T)}^2 + \norm{\partial_t y_{\widehat\TT}-F(y_{\widehat\TT})}{L^2(T)}^2\\
 &\lesssim  \eta^2_\TT(T)+ \eta^2_{\widehat\TT}(T)+\norm{y_\TT-y_{\widehat\TT}}{L^2(T)}^2.
\end{align*}
Note that we slightly abuse notation to write $\eta_{\widehat \TT}(T)^2 \coloneq \sum_{\widehat T \in \widehat \TT: \widehat T \subset T} \eta_{\widehat \TT}(\widehat T)^2$.
\revision{Summing over all $T \in \TT$} using the two preceding estimates yields
\begin{align}\label{eq:h1_esti}
    \begin{split}
        \norm{\partial_t y_\TT-\partial_t y_{\widehat\TT}}{L^2(t_0, \tend)}^2 &\lesssim  \sum_{T \in \TT \setminus \widehat\TT}\eta_{\TT}(T)^2 \\
 &\quad +\sum_{T \in \widehat\TT \setminus \TT}\eta_{\widehat\TT}(T)^2+ \norm{y_\TT-y_{\widehat\TT}}{L^2(t_0, \tend)}^2.
    \end{split}
\end{align}
\revision{To conclude, we need a corresponding upper bound} for $ \norm{y_\TT-y_{\widehat\TT}}{L^2(t_0, \tend)}^2$. \revision{For this, we invoke} Lemma~\ref{lem:help_disc_rel}. \revision{As all the assumptions are satisfied in our setting, we get}
\begin{align*}
        |y_{\widehat\TT}(t)-y_\TT(t)| \lesssim \sum_{T \in \widehat\TT \setminus \TT} \int_T \Big(|R(s)|+|\widehat{R}(s)|\Big)\,\mathrm{d}s.
    \end{align*}
    Squaring and integrating from $t_0$ to $\tend$ and using that \revision{the estimator bounds the residual from above} gives
    \begin{align}\label{eq:l2_estimate1}
        \norm{y_\TT-y_{\widehat\TT}}{L^2(t_0, \tend)}^2 \lesssim \sum_{T \in \TT \setminus \widehat\TT}\eta_{\TT}(T)^2+\sum_{T \in \widehat\TT \setminus \TT}\eta_{\widehat\TT}(T)^2.
    \end{align}
    By inserting \eqref{eq:l2_estimate1} into \eqref{eq:h1_esti} one obtains
    \begin{align*}
    \norm{\partial_t y_\TT-\partial_t y_{\widehat\TT}}{L^2(t_0, \tend)}^2 &\lesssim  \sum_{T \in \TT \setminus \widehat\TT}\eta_{\TT}(T)^2 \\
 &\quad +\sum_{T \in \widehat\TT \setminus \TT}\eta_{\widehat\TT}(T)^2.
\end{align*}
    The strong reduction estimate from Lemma~\ref{lem:red_1} yields
\begin{align*}
    \norm{\partial_t y_\TT-\partial_t y_{\widehat\TT}}{L^2(t_0, \tend)}^2 &\leq  C\sum_{T \in \TT \setminus \widehat\TT}\eta_{\TT}(T)^2 \\
 &+C C_{\rm red}\norm{h_\TT}{L^\infty(t_0,\tend)}^2\norm{y_{\TT}-y_{\widehat\TT}}{H^1(t_0,\tend)}^2,
\end{align*}
for some constant $C>1$. Note that a Poincaré inequality and the fact that $y_\TT(t_0)=y_{\widehat \TT}(t_0)$ give the equivalence of $ \norm{\partial_t (\cdot)}{L^2(t_0, \tend)}$ and $\norm{\cdot}{H^1(t_0,\tend)}$, with constant  $C_{\rm p}$. If the initial mesh is fine enough such that $C_{\rm red} C_{\rm p}\norm{h_\TT}{L^\infty(t_0,\tend)}^2<\frac{1}{C}$, one can move the corresponding term to the \revision{left-hand side} to get
\begin{align*}
    \norm{\partial_t y_\TT-\partial_t y_{\widehat\TT}}{L^2(t_0, \tend)}^2 \lesssim \sum_{T \in \TT \setminus \widehat\TT}\eta_{\TT}(T)^2,
\end{align*}
which concludes the proof. 
\end{proof}
\subsection{\revision{Higher-order} discretizations $p>1$}\noindent\label{sec:high_order}
\revision{The explicit structure of the scheme in \eqref{eq:ode_disc} makes it well-defined} \jn{for arbitrary order $p\ge1$}. However, we will restrict ourselves to the case of affine right-hand sides to be able to perform a rigorous analysis. The analysis for general Lipschitz continuous right-hand sides remains an interesting open problem. \revision{In the following, we consider}
\begin{align*}
    F(y)=Ay+b,
\end{align*}
with $A \in \mathbb{R}^{d \times d}$ and $b \in \mathbb{R}^d$. \revision{Our method requires an explicit lower-order time-stepping scheme, which we construct recursively.} Namely, for $p=2$, we choose $I^F_T$ as the explicit Euler scheme \revision{investigated in the previous section}. This means
\begin{align*}
    I^F_T(z)(t)=z+(t-t_i)F(z), \quad t \in T=[t_i,t_{i+1}].
\end{align*}
As this will lead to an explicit time-stepping scheme, we can reuse it for $p=3$ and continue recursively. \revision{For this simple right-hand side, we can give explicit expressions for the resulting schemes. We collect these in the next lemma.}
\begin{lemma}
    The recursive procedure to obtain \revision{higher-order} explicit \revision{schemes described above yields} the following approximation for $p\geq1$ and $F(y)=Ay+b$:
    \begin{align}
    \label{eq:lin_scheme}
    \begin{split}
    y_\TT(t_0) &= y_0,\\
    y_\TT(t) &= \sum_{k=1}^{p} \frac{1}{k!}A^{k-1}(Ay_\TT(t_i)+b)(t-t_i)^k+ y_\TT(t_i), \quad t \in [t_i,t_{i+1}].
    \end{split}
\end{align}
Furthermore, for $T=[t_i,t_{i+1}]$, the residual  $\big(\partial_t y_\TT-Ay_\TT-b\big)_{|T}$ has a zero of order $p$ at $t_i$, i.e.
    \begin{align}\label{eq:lin_res}
        \partial_t y_\TT(t)-Ay_\TT(t)-b=-\frac{1}{p!} A^p (A y_\TT(t_i)+b)(t-t_i)^p, \quad t \in T.
    \end{align}
\end{lemma}
\begin{proof}
    We will use \revision{induction on} $p\geq1$. The base case $p=1$ is clear by definition. Suppose the assertion is true for some $p\geq 1$. It holds for ${y_\TT}_{|T} \in P^{p+1}(T)$ that
    \begin{align*}
        \int_T \Big( \partial_t y_\TT - F\big(I^F_T(y_\TT(t_i))\big) \Big) \cdot v \,\mathrm{d}t &= 0, \quad v \in P^p(T)
    \end{align*}
    where $I^F_T(y_\TT(t_i))(t)=\sum_{k=1}^{p} \frac{1}{k!}A^{k-1}(Ay_\TT(t_i)+b)(t-t_i)^k+ y_\TT(t_i)$. \revision{The affine structure of $F$ allows us to test} with $v=\partial_t y_\TT - F\big(I^F_T(y_\TT(t_i))\big)$, which yields that
    \begin{align*}
        \partial_t y_\TT =F\big(I^F_T(y_\TT(t_i))\big)=A\sum_{k=1}^{p} \frac{1}{k!}A^{k-1}(Ay_\TT(t_i)+b)(t-t_i)^k+ Ay_\TT(t_i)+b.
    \end{align*}
    Integrating from $t_i$ to $t \in [t_i,t_{i+1}]$ yields
    \begin{align*}
        y_\TT(t)=y_\TT(t_i)+ \sum_{k=1}^{p+1} \frac{1}{k!}A^{k-1}(Ay_\TT(t_i)+b)(t-t_i)^k,
    \end{align*}
    \revision{which proves the first statement. The second statement follows from basic manipulations of \eqref{eq:lin_scheme}, concluding the proof.}
\end{proof}
This allows us to work with an explicit expression for $y_\TT$ in the following. Investigating \eqref{eq:lin_scheme} and \eqref{eq:lin_res} gives the following \revision{elementwise a posteriori error estimator}
\begin{align}\label{eq:lin_esti}
\eta_\TT(T) \coloneq |T|^p \norm{A \partial_t^p y_\TT}{L^2(T)} \simeq
    \norm{R}{L^2(T)}
    .
\end{align}
In the following we aim to reproduce the analysis of the previous section for the \revision{higher-order} schemes presented here. \revision{For this, we need} \textit{strong} stability and reduction estimates. They are stated in \revision{the following lemmas}.
\begin{lemma}\label{lem:stab_2}
      The error estimator defined in \eqref{eq:lin_esti} satisfies \textit{strong} stability in the sense that for a refinement $\widehat\TT$ of $\TT$ it holds
    \begin{align*}
        \Bigg| \Big(\sum_{T \in \TT \cap \widehat\TT}\eta_{\TT}(T)^2 \Big)^{1/2} &-\Big(\sum_{T \in \TT \cap \widehat\TT} \eta_{\widehat\TT}(T)^2\Big)^{1/2} \Bigg| \\
 &\leq  C_{\rm stab} \norm{h_\TT}{L^\infty(t_0,\tend)}\norm{y_{\TT}-y_{\widehat\TT}}{H^1(t_0,\tend)},
    \end{align*}
     where $y_\TT$ and $y_{\widehat{\TT}}$ are discrete solutions on the respective meshes, $h_\TT(t)= |T|$ for $t \in \text{int}(T)$ is the mesh-size function and $C_{\rm stab}>1$.
\end{lemma}
\begin{proof}
    Triangle inequalities show
    \begin{align*}
        \Bigg| \Big(\sum_{T \in \TT \cap \widehat\TT}\eta_{\TT}(T)^2 \Big)^{1/2} &-\Big(\sum_{T \in \TT \cap \widehat\TT} \eta_{\widehat\TT}(T)^2\Big)^{1/2} \Bigg|^2
        \\
 &\leq  \sum_{T \in \TT \cap \widehat\TT}  |T|^{2p} \norm{A\big(\partial_t^p y_\TT-\partial_t^p y_{\widehat\TT}\big)}{L^2(T)}^2\\
 &\lesssim  \norm{A}{}^2 \norm{h_\TT}{L^\infty(t_0,\tend)}^2 \norm{\partial_ty_\TT-\partial_t y_{\widehat\TT}}{L^2(t_0,\tend)}^2,
    \end{align*}
    where we use an inverse inequality in the last step.
\end{proof}
\begin{lemma}\label{lem:red_2}
    The error estimator defined in \eqref{eq:lin_esti} satisfies \textit{strong} reduction in the sense that for a refinement $\widehat\TT$ of $\TT$ it holds
    \begin{align*}
    \sum_{T \in \widehat\TT \setminus \TT} \eta_{\widehat\TT}(T)^2 &\leq  q \sum_{T \in  \TT\setminus \widehat\TT} \eta_{\TT}(T)^2 \\
 &\quad + C_{\rm red}\norm{h_\TT}{L^\infty(t_0,\tend)}^2\norm{y_{\TT}-y_{\widehat\TT}}{H^1(t_0,\tend)}^2,
\end{align*}
    for $0<q<1$ and $C_{\rm red}>1$.
\end{lemma}
\begin{proof}
    Let $T \in \TT \setminus \widehat\TT$ and  $\{\widehat{T}_i\}_{i=1}^n \subset \widehat\TT \setminus \TT$ such that $T=\bigcup_{i=1}^n \widehat{T}_i $. It holds
    \begin{align*}
       \sum_{i=1}^n \eta_{\widehat\TT}(\widehat{T}_i)^2 &=  \sum_{i=1}^n |\widehat{T}_i|^{2p}\norm{A\partial_t^p y_{\widehat\TT}}{L^2(\widehat{T}_i)}^2\\
 &\leq  \frac{|T|^{2p}}{4}\norm{A\partial_t^p y_{\TT}}{L^2(T)}^2 + C_{\rm inv}\frac{|T|^{2}}{4}\norm{A\partial_t y_{\TT}-A\partial_t y_{\widehat\TT}}{L^2(T)}^2 \\
 &\leq  \frac14 \eta_\TT(T)^2 + C_{\rm inv}\norm{A}{}^2\frac{|T|^{2}}{4}\norm{\partial_t y_{\TT}-\partial_t y_{\widehat\TT}}{L^2(T)}^2,
    \end{align*}
    \revision{since} $|\widehat{T}_i|\leq \frac12 |T|$ and by an inverse inequality in the second step. \revision{Summing over all $T \in \TT \setminus \widehat\TT$} yields the result.
\end{proof}
As before, linear convergence of the estimator follows immediately from Lemmas~\ref{lem:stab_2} and~\ref{lem:red_2}. We will state the result in the following. Note that the proof is exactly the same as \revision{in the lowest-order case}, so we omit it.
\begin{lemma}
    \label{lem:rlinearconv_lin}
    Let $\TT_0, \TT_{1},\ldots$ be meshes generated by Algorithm~\ref{alg:adaptive} with \revision{corresponding discrete solutions} $y_{0}, y_{1},\ldots$ and estimators $\eta_{0}, \eta_{1},\ldots$. Suppose the initial discretization $\TT_0$ is fine enough. \revision{Furthermore, suppose} that the estimator satisfies strong stability and reduction as outlined in Lemmas~\ref{lem:stab_2} and~\ref{lem:red_2}. Then there is a $0<\kappa<1$ \revision{such that, for all $k,\ell \in \N$,}
    \begin{align*}
        \eta^2_{\ell+k} \leq \kappa^k \eta^2_\ell.
    \end{align*}
\end{lemma}
The only link missing in order to prove optimal convergence rates is the discrete reliability estimate, already mentioned above. The proof \revision{follows the same lines as the proof} of Lemma~\ref{lem:discrel}.
\begin{lemma}\label{lem:discrel2}
    Let $\widehat\TT$ be a refinement of $\TT$. Under the assumption of a fine enough initial mesh $\TT_0$, \revision{discrete reliability holds in the sense}
        \begin{align*}
        \norm{y_{\widehat\TT}-y_{\TT}}{H^1(t_0,\tend)}^2 \leq C_{\rm drel}^2 \sum_{T \in \TT \setminus \widehat\TT}\eta_{\TT}(T)^2,
    \end{align*}
    for corresponding discrete solutions $y_{\widehat\TT}, y_\TT$ defined in \eqref{eq:lin_scheme} and $C_{\rm drel}>1.$
\end{lemma}
\begin{proof}
 We first consider $[t_i, t_{i+1}]=T \in \TT \cap \widehat\TT$. From the definition of the numerical scheme \eqref{eq:lin_scheme} we get that
    \begin{align*}
        |\partial_t y_\TT-\partial_t y_{\widehat\TT}|=\Big|\sum_{k=1}^{p} \frac{1}{(k-1)!}A^{k}(y_\TT(t_i)-y_{\widehat\TT}(t_i))(t-t_i)^{k-1}\Big|,
    \end{align*}
    which yields \dn{by using norm equivalence on \revision{finite-dimensional spaces} followed by scaling\revision{} that}
    \begin{align*}
    \norm{\partial_t y_\TT-\partial_t y_{\widehat\TT}}{L^2(T)} \lesssim \norm{y_\TT-y_{\widehat\TT}}{L^2(T)},
\end{align*}
where the implied constant depends on $\norm{A}{}$ and $p$.
Furthermore, we also get the estimate \eqref{eq:initial_esti} by using similar arguments. For $[t_i,t_{i+1}]=T \in \TT \setminus \widehat\TT$ we get
\begin{align*}
 &\norm{\partial_t y_\TT-\partial_t y_{\widehat\TT}}{L^2(T)}^2\ \\
 &\lesssim  \norm{\partial_t y_\TT-F(y_\TT)}{L^2(T)}^2+\norm{F(y_\TT)-F(y_{\widehat\TT})}{L^2(T)}^2 + \norm{\partial_t y_{\widehat\TT}-F(y_{\widehat\TT})}{L^2(T)}^2\\
 &\lesssim  \eta^2_\TT(T)+ \eta^2_{\widehat\TT}(T)+\norm{y_\TT-y_{\widehat\TT}}{L^2(T)}^2.
\end{align*}
\revision{Summing over all $T \in \TT$} using the two preceding estimates yields
\begin{align}\label{eq:h1_esti_lin}
    \begin{split}
        \norm{\partial_t y_\TT-\partial_t y_{\widehat\TT}}{L^2(t_0, \tend)}^2 &\lesssim  \sum_{T \in \TT \setminus \widehat\TT}\eta_{\TT}(T)^2 \\
 &\quad +\sum_{T \in \widehat\TT \setminus \TT}\eta_{\widehat\TT}(T)^2+ \norm{y_\TT-y_{\widehat\TT}}{L^2(t_0, \tend)}^2.
    \end{split}
\end{align}
\revision{To bound the term $\norm{y_\TT-y_{\widehat\TT}}{L^2(t_0, \tend)}^2$,} we invoke Lemma~\ref{lem:help_disc_rel}. \revision{As all the assumptions are satisfied in our setting, we get}
\begin{align*}
        |y_{\widehat\TT}(t)-y_\TT(t)| \lesssim \sum_{T \in \widehat\TT \setminus \TT} \int_T \Big(|R(s)|+|\widehat{R}(s)|\Big)\,\mathrm{d}s.
    \end{align*}
    Squaring and integrating from $t_0$ to $\tend$ and using that \revision{the estimator bounds the residual from above} gives
    \begin{align*}
        \norm{y_\TT-y_{\widehat\TT}}{L^2(t_0, \tend)}^2 \lesssim \sum_{T \in \TT \setminus \widehat\TT}\eta_{\TT}(T)^2+\sum_{T \in \widehat\TT \setminus \TT}\eta_{\widehat\TT}(T)^2.
    \end{align*}
    \revision{Inserting this into \eqref{eq:h1_esti_lin}} we obtain
    \begin{align*}
    \norm{\partial_t y_\TT-\partial_t y_{\widehat\TT}}{L^2(t_0, \tend)}^2 &\lesssim  \sum_{T \in \TT \setminus \widehat\TT}\eta_{\TT}(T)^2 \\
 &\quad +\sum_{T \in \widehat\TT \setminus \TT}\eta_{\widehat\TT}(T)^2.
\end{align*}
    The strong reduction estimate from Lemma~\ref{lem:red_2} yields
\begin{align*}
 &\norm{\partial_t y_\TT-\partial_t y_{\widehat\TT}}{L^2(t_0, \tend)}^2 \\
 &\leq  C \sum_{T \in \TT \setminus \widehat\TT}\eta_{\TT}(T)^2 + C C_{\rm red}\norm{h_\TT}{L^\infty(t_0,\tend)}^2\norm{y_{\TT}-y_{\widehat\TT}}{H^1(t_0,\tend)}^2,
\end{align*}
for a constant $C>1$.
 Note that a Poincaré inequality and the fact that $y_\TT(t_0)=y_{\widehat \TT}(t_0)$ give the equivalence of $ \norm{\partial_t (\cdot)}{L^2(t_0, \tend)}$ and $\norm{\cdot}{H^1(t_0,\tend)}$ with constant  $C_{\rm p}$. If the initial mesh is fine enough such that $C_{\rm red} C_{\rm p}\norm{h_\TT}{L^\infty(t_0,\tend)}^2<\frac{1}{C}$, one can move the corresponding term to the \revision{left-hand side} to get
\begin{align*}
    \norm{\partial_t y_\TT-\partial_t y_{\widehat\TT}}{L^2(t_0, \tend)}^2 \lesssim \sum_{T \in \TT \setminus \widehat\TT}\eta_{\TT}(T)^2,
\end{align*}
which concludes the proof.
\end{proof}
All the tools needed to prove Theorem~\ref{thm:optim} have been collected. \revision{As the analysis is similar for} \revision{the lowest-order nonlinear and higher-order linear cases}, \revision{we treat both cases simultaneously}.
\begin{proof}[Proof of Theorem~\ref{thm:optim}]
Let $0<\theta<1$ be given.
    \revision{We use Lemmas~\ref{lem:stab_1} and~\ref{lem:discrel} or Lemmas~\ref{lem:stab_2} and~\ref{lem:discrel2}, respectively} to apply \cite[Proposition 4.12]{axioms}. We choose $\TT_0$ sufficiently fine such that $\theta<\theta_\star:=(1+\norm{h_{\TT_0}}{L^\infty(t_0,\tend)}^2C_{\rm stab}^2 C^2_{\rm drel})^{-1}$ and such that \revision{all the lemmas} are applicable. \revision{This allows us to further apply} \cite[Lemma 4.14]{axioms}. Due to Lemma~\ref{lem:rlinearconv} or Lemma~\ref{lem:rlinearconv_lin} we can apply  \cite[Proposition 4.15]{axioms} with $\theta_0= \theta$ which yields the result and concludes the proof.
\end{proof}
\section{A posteriori error control for the $W^{1,\infty}$-error}
\revision{In this section, we discuss how the estimators obtained above can also be used to control the error in the $W^{1,\infty}$-norm.} The discussion is similar to the previous one. \revision{We nevertheless point out} the major differences that occur in the analysis and the algorithm. \revision{Lemma~\ref{lem:rel} already shows that the $L^\infty$-norm of the residual provides an upper bound for} $\norm{y-y_\TT}{W^{1,\infty}(t_0,\tend)}$. The different norms naturally give rise to \revision{residual-based} elementwise error estimators of the form
\begin{align*}
    \eta_\TT^\infty=\max_{T \in \TT} \eta_\TT^\infty(T),
\end{align*}
\revision{which lead to an adaptive algorithm using the maximum strategy for refinement in place of Dörfler marking.}
\begin{algorithm}[h]
    \caption{Adaptive Algorithm 2}
    \label{alg:adaptivemax}
    \begin{algorithmic}[1]
       \item \textbf{Input:} Initial mesh $\TT_0$, parameter $0<\theta<1$.
        \For{ $\ell=0,1,2,\ldots$}:
        \State Compute $y_{\TT_\ell}$ from \eqref{eq:ode_disc}.
        \State Compute error estimates $\eta_{\TT_\ell}^\infty(T)$ for all $T \in \TT_\ell$.
        \If{$\eta^\infty_{\TT_\ell}=0$} \State break. \EndIf
        \State Mark all $T \in \TT_\ell$ for refinement that satisfy
        \begin{align*}
            \eta_{\TT_\ell}^\infty(T)\geq (1-\theta) \eta_{\TT_\ell}^\infty.
        \end{align*}
        \State Bisect all marked elements to obtain a new mesh $\TT_{\ell+1}$.
        \EndFor
    \end{algorithmic}
\end{algorithm}
\begin{theorem}\label{thm:optim_max}
\revision{Assume the setting of Section~}\ref{sec:low_order_inf} or~\ref{sec:high_order_inf}.
    Let $0<\theta<1$ be arbitrary but given. Suppose the initial discretization $\TT_0$ is fine enough. Let $s>0$ be such that
    \begin{align*}
        C_{\rm best}:=\sup_{N \in \mathbb{N}} \inf_{\substack{\TT \in \T \\ \#\TT-\#\TT_0\leq N}}(N+1)^s \eta_\TT^\infty<\infty.
    \end{align*}
    The sequence of meshes $\TT_\ell$ obtained from Algorithm~\ref{alg:adaptivemax}
    satisfies quasi-optimal convergence rates in the sense
    \begin{align}
        \sup_{\ell \in \mathbb{N}} \eta_{\TT_\ell}^\infty(\#\TT_\ell-\#\TT_0+1)^s \leq C_{\rm opt} C_{\rm best}.
    \end{align}
\end{theorem}
\subsection{\revision{Lowest-order} discretization $p=1$}\label{sec:low_order_inf}
\revision{The approximation scheme for $p=1$ remains unchanged, i.e.,} the explicit Euler method is used. The previous discussion allows us to derive an \revision{error estimator} that fits in this maximum setting. It holds for every $T \in \TT$ that
\begin{align}\label{eq:esti_inf}
\eta_\TT^\infty(T) \coloneqq L|T|\norm{\partial_t y_\TT}{L^\infty(T)}=
L |T| |F(y_\TT(t_i))| \geq
    \norm{R}{L^\infty(T)} .
\end{align}
The \textit{strong} stability and reduction estimates \revision{follow immediately from the same strategy} as in the previous section. \revision{As they have a slightly different structure, we state} them in a unified manner in the following.
\begin{lemma}\label{lem:stabred_max}
    The error estimator defined in \eqref{eq:esti_inf} satisfies strong stability and reduction, in the sense that for a refinement $\widehat\TT$ of $\TT$ it holds that
    \begin{align*}
        \Big| \max_{T \in \TT \cap \widehat\TT} \eta_\TT^\infty(T)- \max_{T \in \TT \cap \widehat\TT} \eta_{\widehat\TT}^\infty(T) \Big| \leq C_{\rm stab} \norm{h_\TT}{L^\infty(t_0,\tend)}\norm{y_{\TT}-y_{\widehat\TT}}{W^{1,\infty}(t_0,\tend)},
    \end{align*}
    and
    \begin{align*}
    \max_{T \in \widehat\TT \setminus \TT} \eta_{\widehat\TT}^\infty(T) &\leq  q  \max_{T \in \TT \setminus \widehat\TT} \eta_{\TT}^\infty(T) \\
 &\quad + C_{\rm red} \norm{h_\TT}{L^\infty(t_0,\tend)}\norm{y_{\TT}-y_{\widehat\TT}}{W^{1,\infty}(t_0,\tend)},
\end{align*}
    where $0<q<1$ and $C_{\rm red},C_{\rm stab}>1$.
\end{lemma}
Linear convergence can be obtained from Lemma~\ref{lem:stabred_max} \revision{as before} and is stated \revision{in the next lemma}.
\begin{lemma}\label{lem:linconv_inf}
    Let $\TT_0, \TT_{1},\ldots$ be meshes generated by Algorithm~\ref{alg:adaptivemax} with \revision{corresponding discrete solutions} $y_{0}, y_{1},\ldots$ and estimators $\eta_{0}^\infty, \eta_{1}^\infty,\ldots$. Suppose the initial discretization $\TT_0$ is fine enough. \revision{Furthermore, suppose} that the estimator satisfies strong stability and reduction as outlined in Lemma~\ref{lem:stabred_max}. Then there is a $0<\kappa<1$ \revision{such that, for all $k,\ell \in \N$,}
    \begin{align*}
        \eta_{\ell+k}^\infty \leq \kappa^k \eta_\ell^\infty.
    \end{align*}
\end{lemma}
\begin{proof}
     Let $\ell \in \N$. By using stability and reduction, \revision{i.e.,} Lemma~\ref{lem:stabred_max}, we obtain
        \begin{align*}
            \eta_{\ell+1}^\infty &=  \max\Big\{\max_{T \in \TT_{\ell+1} \setminus \TT_{\ell}} \eta_{\ell+1}^\infty(T), \max_{T \in \TT_{\ell+1} \cap \TT_{\ell}} \eta_{\ell+1}^\infty(T)\Big\}\\
 &\leq  \max \Big\{q \max_{T \in \TT_{\ell} \setminus \TT_{\ell+1}} \eta_{\ell}^\infty(T),  \max_{T \in \TT_{\ell+1} \cap \TT_{\ell}} \eta_{\ell}^\infty(T)\Big\}\\
 &+ C_{\rm est} \norm{h_{\TT_\ell}}{L^\infty(t_0,\tend)}\norm{y_{\ell}-y_{\ell+1}}{W^{1,\infty}(t_0,\tend)}
        \end{align*}
        with $C_{\rm est}=\max\{C_{\rm red},C_{\rm stab}\}$. Using the marking strategy with $\theta \in (0,1)$ yields
        \begin{align*}
             \eta_{\ell+1}^\infty \leq \max\{q,1-\theta\}\eta^\infty_\ell + C_{\rm est} \norm{h_{\TT_\ell}}{L^\infty(t_0,\tend)}\norm{y_{\ell}-y_{\ell+1}}{W^{1,\infty}(t_0,\tend)}.
        \end{align*}
        We have that $q_{\rm red}:= \max\{q,1-\theta\}<1$. Since we have the reliability estimate
        \begin{align*}
            \norm{y_{\ell}-y_{\ell+1}}{W^{1,\infty}(t_0,\tend)} \leq C_{\rm rel} \Big(  \eta^\infty_\ell+\eta^\infty_{\ell+1} \Big),
        \end{align*}
        we get for $C=C_{\rm rel} C_{\rm est}$ that
        \begin{align*}
            \eta_{\ell+1}^\infty (1-C\norm{h_{\TT_\ell}}{L^\infty(t_0,\tend)}) \leq \eta_{\ell}^\infty (q_{\rm red}+C\norm{h_{\TT_\ell}}{L^\infty(t_0,\tend)}),
        \end{align*}
        which yields under the assumption of a fine initial mesh such that $C \norm{h_{\TT_\ell}}{L^\infty(t_0,\tend)} \leq C \norm{h_{\TT_0}}{L^\infty(t_0,\tend)}<1$ that
        \begin{align*}
            \eta_{\ell+1}^\infty \leq \frac{q_{\rm red}+C\norm{h_{\TT_0}}{L^\infty(t_0,\tend)}}{1-C\norm{h_{\TT_0}}{L^\infty(t_0,\tend)}} \eta_{\ell}^\infty.
        \end{align*}
        Choosing an even finer initial discretization $\TT_0$ such that
        \begin{align*}
            \kappa \coloneq \frac{q_{\rm red}+C\norm{h_{\TT_0}}{L^\infty(t_0,\tend)}}{1-C\norm{h_{\TT_0}}{L^\infty(t_0,\tend)}}<1,
        \end{align*}
        concludes the proof.
\end{proof}
\revision{The discrete reliability estimate can be derived} \revision{in a manner similar to that in} the preceding section by using Lemma~\ref{lem:help_disc_rel} and the strong reduction estimate.
\begin{lemma}\label{lem:discrel_inf}
    Let $\widehat\TT$ be a refinement of $\TT$. Under the assumption of a fine enough initial mesh $\TT_0$, \revision{discrete reliability holds in the sense}
        \begin{align*}
        \norm{y_{\widehat\TT}-y_{\TT}}{W^{1,\infty}(t_0,\tend)} \leq C_{\rm drel} \max_{T \in \TT \setminus \widehat\TT}\eta_{\TT}^\infty(T),
    \end{align*}
    for corresponding discrete solutions $y_{\widehat\TT}, y_\TT$ defined in \eqref{eq:exp_euler} and $C_{\rm drel}>1.$
\end{lemma}
\begin{proof}
 We first consider $[t_i, t_{i+1}]=T \in \TT \cap \widehat\TT$. From the definition of the numerical scheme \eqref{eq:exp_euler} we get that
    \begin{align*}
        |\partial_t y_\TT-\partial_t y_{\widehat\TT}|=\Big|F(y_\TT(t_i))-F(y_{\widehat\TT}(t_i))\Big|,
    \end{align*}
    which yields by using the Lipschitz condition on $F$ that
    \begin{align*}
    \norm{\partial_t y_\TT-\partial_t y_{\widehat\TT}}{L^\infty(T)} \lesssim \norm{y_\TT-y_{\widehat\TT}}{L^\infty(T)}.
\end{align*}
Furthermore, we also get the estimate \eqref{eq:initial_esti} by using similar arguments. For $[t_i,t_{i+1}]=T \in \TT \setminus \widehat\TT$ we get
\begin{align*}
 &\norm{\partial_t y_\TT-\partial_t y_{\widehat\TT}}{L^\infty(T)}\ \\
 &\leq  \norm{\partial_t y_\TT-F(y_\TT)}{L^\infty(T)}+\norm{F(y_\TT)-F(y_{\widehat\TT})}{L^\infty(T)} + \norm{\partial_t y_{\widehat\TT}-F(y_{\widehat\TT})}{L^\infty(T)}\\
 &\lesssim  \eta^\infty_\TT(T)+ \eta^\infty_{\widehat\TT}(T)+\norm{y_\TT-y_{\widehat\TT}}{L^\infty(T)},
\end{align*}
where $\eta^\infty_{\widehat\TT}(T)=\max_{\widehat{T}\subset T}\eta^\infty_{\widehat\TT}(\widehat{T}).$
Maximizing over all $T \in \TT$ using the two preceding estimates yields
\begin{align}\label{eq:h1_esti_inf}
    \begin{split}
        \norm{\partial_t y_\TT-\partial_t y_{\widehat\TT}}{L^\infty(t_0, \tend)} &\lesssim  \max_{T \in \TT \setminus \widehat\TT}\eta_{\TT}^\infty(T) \\
 &\quad +\max_{T \in \widehat\TT \setminus \TT}\eta_{\widehat\TT}^\infty(T)+ \norm{y_\TT-y_{\widehat\TT}}{L^\infty(t_0, \tend)}.
    \end{split}
\end{align}
To estimate the $L^\infty(t_0,\tend)$ norm of $y_\TT-y_{\widehat \TT}$ we invoke Lemma~\ref{lem:help_disc_rel}. \revision{As all the assumptions are satisfied in our setting, we get}
\begin{align*}
        |y_{\widehat\TT}(t)-y_\TT(t)| \lesssim \sum_{T \in \widehat\TT \setminus \TT} \int_T \Big(|R(s)|+|\widehat{R}(s)|\Big)\,\mathrm{d}s.
    \end{align*}
    \revision{Taking the supremum over the domain and using that} \revision{the estimator bounds the residual from above} gives
    \begin{align*}
    \norm{y_\TT-y_{\widehat\TT}}{L^\infty(t_0, \tend)} &\lesssim  \max_{T \in \TT \setminus \widehat\TT}\eta_{\TT}^\infty(T) \\
 &\quad +\max_{T \in \widehat\TT \setminus \TT}\eta_{\widehat\TT}^\infty(T).
\end{align*}
    \revision{Inserting this into \eqref{eq:h1_esti_inf}} we obtain
    \begin{align*}
    \norm{\partial_t y_\TT-\partial_t y_{\widehat\TT}}{L^\infty(t_0, \tend)} &\lesssim  \max_{T \in \TT \setminus \widehat\TT}\eta_{\TT}^\infty(T) \\
 &\quad +\max_{T \in \widehat\TT \setminus \TT}\eta_{\widehat\TT}^\infty(T).
\end{align*}
    The strong reduction estimate from Lemma~\ref{lem:stabred_max} yields
\begin{align*}
 &\norm{\partial_t y_\TT-\partial_t y_{\widehat\TT}}{L^\infty(t_0, \tend)} \\
 &\leq  C \max_{T \in \TT \setminus \widehat\TT}\eta_{\TT}^\infty(T) + C C_{\rm red}\norm{h_\TT}{L^\infty(t_0,\tend)}\norm{y_{\TT}-y_{\widehat\TT}}{W^{1,\infty}(t_0,\tend)},
\end{align*}
for $C>1$.
Note that a Poincaré inequality and the fact that $y_\TT(t_0)=y_{\widehat \TT}(t_0)$ give the equivalence of $ \norm{\partial_t (\cdot)}{L^\infty(t_0, \tend)}$ and $\norm{\cdot}{W^{1,\infty}(t_0,\tend)}$, with constant  $C_{\rm p}$. If the initial mesh is fine enough such that $C_{\rm red} C_{\rm p}\norm{h_\TT}{L^\infty(t_0,\tend)}<\frac{1}{C}$, one can move the corresponding term to the \revision{left-hand side} to get
\begin{align*}
    \norm{ y_\TT- y_{\widehat\TT}}{W^{1,\infty}(t_0, \tend)} \lesssim \max_{T \in \TT \setminus \widehat\TT}\eta_{\TT}^\infty(T),
\end{align*}
which concludes the proof.
\end{proof}
\subsection{\revision{Higher-order} discretizations $p>1$} \label{sec:high_order_inf}
\jn{The $W^{1,\infty}$ analysis of the lowest-order case carries over to the
higher-order scheme \eqref{eq:lin_scheme} with affine right-hand side
$F(y)=Ay+b$ with only minor modifications.} Therefore, we just state the results needed to prove optimal convergence rates as laid out in Theorem~\ref{thm:optim_max}. 
\jn{We start with the definition of the error estimator and recall the approximation $y_\TT$}
\begin{align}
     y_\TT(t) &= \sum_{k=1}^{p} \frac{1}{k!}A^{k-1}(Ay_\TT(t_i)+b)(t-t_i)^k+ y_\TT(t_i), \quad t \in T\\
    \label{eq:lin_esti_inf}
    \eta_\TT^\infty(T) &\coloneq |T|^p \norm{A \partial_t^p y_\TT}{L^\infty(T)} \simeq
    \norm{R}{L^\infty(T)}
\end{align}
\begin{lemma}\label{lem:stabredmax_lin}
The error estimator defined in \eqref{eq:lin_esti_inf} satisfies strong stability and reduction, in the sense that for a refinement $\widehat\TT$ of $\TT$ it holds that
    \begin{align*}
        \Big| \max_{T \in \TT \cap \widehat\TT} \eta_\TT^\infty(T)- \max_{T \in \TT \cap \widehat\TT} \eta_{\widehat\TT}^\infty(T) \Big| \leq C_{\rm stab} \norm{h_\TT}{L^\infty(t_0,\tend)}\norm{y_{\TT}-y_{\widehat\TT}}{W^{1,\infty}(t_0,\tend)},
    \end{align*}
    and
    \begin{align*}
    \max_{T \in \widehat\TT \setminus \TT} \eta_{\widehat\TT}^\infty(T) &\leq  q  \max_{T \in \TT \setminus \widehat\TT} \eta_{\TT}^\infty(T) \\
 &\quad + C_{\rm red} \norm{h_\TT}{L^\infty(t_0,\tend)}\norm{y_{\TT}-y_{\widehat\TT}}{W^{1,\infty}(t_0,\tend)},
\end{align*}
    where $0<q<1$ and $C_{\rm red},C_{\rm stab}>1$.
    \revision{Assume additionally that} the initial discretization $\TT_0$ is fine enough. Let $\widehat\TT$ be a refinement of $\TT$. \revision{Discrete reliability holds in the sense}
        \begin{align*}
        \norm{y_{\widehat\TT}-y_{\TT}}{W^{1,\infty}(t_0,\tend)} \leq C_{\rm drel} \max_{T \in \TT \setminus \widehat\TT}\eta_{\TT}^\infty(T),
    \end{align*}
    for corresponding discrete solutions $y_{\widehat\TT}, y_\TT$ defined in \eqref{eq:lin_scheme} and $C_{\rm drel}>1.$
    Denote by $\TT_0, \TT_{1},\ldots$ meshes generated by Algorithm~\ref{alg:adaptivemax} with \revision{corresponding discrete solutions} $y_{0}, y_{1},\ldots$ and estimators $\eta_{0}^\infty, \eta_{1}^\infty,\ldots$. \revision{There is a constant $0<\kappa<1$} \revision{such that, for all $k,\ell \in \N$,}
    \begin{align*}
        \eta_{\ell+k}^\infty \leq \kappa^k \eta_\ell^\infty.
    \end{align*}
\end{lemma}
\subsection{Proof of Theorem~\ref{thm:optim_max}}
The aim is to use the framework established in \cite{axioms} and \revision{adapt it wherever necessary to fit our setting}. Linear and nonlinear settings are handled simultaneously as before, \revision{by using the corresponding lemmas in each case}. One has to adapt the \revision{so-called} \textit{optimality of \revision{Dörfler marking}, i.e.,} \cite[Proposition 4.12]{axioms} to the present case of the maximum strategy. That is why we will call the next result \textit{optimality of the maximum strategy}.
\begin{proposition}\label{prop:opti_max} We have
    \begin{enumerate}
        \item[(i)] For all $0<\kappa_0<1$, there exists an initial mesh $\TT_0$ fine enough and a constant $0<\theta_0<1$ such that for all $0<\theta<\theta_0$ and all refinements $\widehat\TT \in \T$ of $\TT \in \T$ it holds if $\eta_\TT^\infty>0$ that
        \begin{align}\label{eq:implication}
            \eta_{\widehat\TT}^\infty \leq \kappa_0 \eta_\TT^\infty \implies (1-\theta) \eta_\TT^\infty > \max_{T \in \TT \cap \widehat\TT} \eta_\TT^\infty(T).
        \end{align}
        \item[(ii)] Suppose $\eta_\TT^\infty>0$. For all $\theta_0$ such that $\theta_\star:=1-\norm{h_{\TT_0}}{L^\infty(t_0,\tend)}C_{\rm stab}C_{\rm drel}>\theta_0>0$, there exists some $0< \kappa_0<1$ such that \eqref{eq:implication} holds for all $0<\theta\leq \theta_0<1$ and all refinements $\widehat\TT \in \T$ of $\TT \in \T$.
    \end{enumerate}
\end{proposition}
\begin{proof}
    (i): It holds by the strong stability estimate Lemma~\ref{lem:stabred_max} or Lemma~\ref{lem:stabredmax_lin} that
    \begin{align*}
    \max_{T \in \TT \cap \widehat\TT} \eta_\TT^\infty(T) &\leq  \max_{T \in \TT \cap \widehat\TT} \eta_{\widehat\TT}^\infty(T) \\
 &\quad + \norm{h_{\TT_0}}{L^\infty(t_0,\tend)}C_{\rm stab} \norm{y_\TT-y_{\widehat\TT}}{W^{1,\infty}(t_0,\tend)},
\end{align*}
    which gives by using the assumption $\eta_{\widehat\TT}^\infty \leq \kappa_0 \eta_\TT^\infty$ and discrete reliability in Lemma~\ref{lem:discrel_inf},
    \begin{align*}
         \max_{T \in \TT \cap \widehat\TT} \eta_\TT^\infty(T) \leq \Big(\kappa_0 +\norm{h_{\TT_0}}{L^\infty(t_0,\tend)}C_{\rm stab}C_{\rm drel}\Big) \eta_\TT^\infty.
    \end{align*}
    Choosing $\theta_0=1-\kappa_0 -\norm{h_{\TT_0}}{L^\infty(t_0,\tend)}C_{\rm stab}C_{\rm drel}>0$, with $\norm{h_{\TT_0}}{L^\infty(t_0,\tend)}$ small enough, which can be controlled by choosing a fine enough initial mesh, yields the result.

    \noindent(ii): Let $1-\norm{h_{\TT_0}}{L^\infty(t_0,\tend)}C_{\rm stab}C_{\rm drel}>\theta_0>0$. If one chooses the reduction factor $\kappa_0<1-\theta_0-\norm{h_{\TT_0}}{L^\infty(t_0,\tend)}C_{\rm stab}C_{\rm drel}$, the result follows from the arguments of (i).
\end{proof} \noindent
\revision{Using Proposition~\ref{prop:opti_max}, one can adapt the crucial results} \cite[Lemma 4.14]{axioms} and \cite[Proposition 4.15]{axioms} by using our Proposition~\ref{prop:opti_max} wherever \cite[Proposition 4.12]{axioms} is used and substituting \textit{maximum strategy} for \textit{Dörfler marking}. \revision{This allows us to prove} Theorem~\ref{thm:optim_max} in the following.
\begin{proof}[Proof of Theorem~\ref{thm:optim_max}]
    Let $0<\theta<1$ be given.
    We use Proposition~\ref{prop:opti_max} and choose $\TT_0$ sufficiently fine such that $0<\theta<\theta_\star:=1-\norm{h_{\TT_0}}{L^\infty(t_0,\tend)}C_{\rm stab}C_{\rm drel}$ and such that \revision{all the lemmas} are applicable. Combining stability, reduction and discrete reliability in a standard fashion yields quasi-monotonicity of the error estimator, \revision{i.e.,} $\eta^\infty_{\widehat \TT} \lesssim \eta^\infty_{\TT}$, for all refinements $\widehat \TT$ of $\TT$. \revision{This allows us to apply} \cite[Lemma 4.14]{axioms}. Due to Lemma~\ref{lem:linconv_inf} or Lemma~\ref{lem:stabredmax_lin} we can apply  \cite[Proposition 4.15]{axioms} with $\theta_0=\theta$ which yields the result and concludes the proof.
\end{proof}
\section{Remarks on the practical implementation}\label{sec:practical}
In this section, we aim to \revision{explain how} the explicit schemes introduced in \eqref{eq:ode_disc} can be interpreted in terms of classical explicit \revision{Runge--Kutta} methods, and thus be implemented as such. {\color{black} The formulation in \eqref{eq:ode_disc} accommodates a broad class of explicit time-stepping schemes through different choices of the local approximation operator. Following the preceding sections, we focus on schemes constructed recursively from the explicit Euler method, which represent only one of many possible constructions within this framework.} \jn{For $p=1$, we retrieve the explicit Euler method.} {\color{black} We continue recursively as described in Section~\ref{sec:high_order}.} \revision{The second-order scheme, i.e., $p=2$, builds} on the explicit Euler method. As the integrals that need to be evaluated to determine $y_\TT$ cannot be computed by hand in general, \revision{one needs to use a quadrature method of sufficiently high order}. \revision{Using the 2-point Radau quadrature with a fixed right endpoint} gives the following conditions for $T=[t_i,t_{i+1}] \in \TT$.
\begin{align*}
  \partial_t y_\TT(t_i+|T|/3) &= F(I^F_T(y_\TT(t_i))(t_i+|T|/3))=F(y_\TT(t_i)+\frac{|T|}{3}F(y_\TT(t_i))):=k_1\\
    \partial_t y_\TT(t_{i+1}) &= F(I^F_T(y_\TT(t_i))(t_{i+1}))=F(y_\TT(t_i)+|T|F(y_\TT(t_i))):=k_2
\end{align*}
We denote the resulting function \revision{again by $y_\TT$}, despite the additional approximation of the integral.
This gives the following expression \revision{for $(y_\TT)_{|T}$ using} the $2$-point Lagrange interpolants $l_{|T|/3} \in P^1(T)$ and $l_{|T|} \in P^1(T)$, 
\begin{align*}
    y_\TT(t)=y_\TT(t_i)+k_1\int_{t_i}^t l_{|T|/3}(s)\,\mathrm{d}s+ k_2\int_{t_i}^t l_{|T|}(s)\,\mathrm{d}s.
\end{align*}
Thus $y_\TT(t_{i+1})$ can be written as a linear combination of \jn{three} stages, resulting in an explicit \revision{Runge--Kutta} scheme, shown in Figure~\ref{fig:rk}. The evaluation at certain intermediate points in order to evaluate the error estimator can also be written as a (different) linear combination of the given stages $k_1$ and $k_2$. Continuing this discussion for $p=3$ and using the 3-point Radau quadrature yields conditions on $\partial_t y_\TT$ at those three points. Those stages only depend on evaluations of the \revision{lower-order} discretization $I_T^F(y_\TT(t_i))$ and thus can be written as linear combinations of the original stages $k_1$ and $k_2$. 
This yields another explicit Runge--Kutta scheme. \revision{The resulting schemes are determined by their Butcher tableaux}, which are shown in Figure~\ref{fig:rk}.
Checking the order conditions yields the following proposition.
\begin{figure}[h]
    \centering
    \begin{subfigure}{0.45\textwidth}
        \centering
        $\begin{array}
{c|c}
0\\
\hline
& 1 
\end{array}$
\caption{$p=1$}
    \end{subfigure}
     \begin{subfigure}{0.45\textwidth}
        \centering
        $\begin{array}
{c|ccc}
0\\
\frac{1}{3} & \frac{1}{3}\\
1 &1 & \\
\hline
& 0 &\frac{3}{4} &\frac{1}{4}
\end{array}$
\caption{$p=2$}
    \end{subfigure}
     \begin{subfigure}{0.9\textwidth}
        \centering
        $\begin{array}
{c|cccccc}
0\\
\frac{1}{3} & \frac{1}{3}\\
1 &1 & \\
\frac{4-\sqrt{6}}{10} & 0 & 0.435-0.09\sqrt{6} & -0.035-0.01 \sqrt{6}\\
\frac{4+\sqrt{6}}{10} & 0 & 0.435+0.09\sqrt{6} & -0.035+0.01 \sqrt{6} \\
1 & 0 & \frac{3}{4} & \frac{1}{4}\\
\hline
& 0 &0 &0 &\frac{4}{9}-\frac{\sqrt{6}}{36} & \frac{4}{9}+\frac{\sqrt{6}}{36}&\frac{1}{9}
\end{array}$
\caption{$p=3$}
    \end{subfigure}
    \caption{Explicit \revision{Runge--Kutta} methods obtained by discretizing the scheme defined in \eqref{eq:ode_disc}.}
    \label{fig:rk}
\end{figure}
\begin{proposition}
    \revision{The explicit RK schemes in Figure~\ref{fig:rk} have convergence orders $1$, $2$, and $3$, respectively.}
\end{proposition}
\revision{Using different quadrature schemes with $p$ points yields different Butcher tableaux} and therefore different explicit schemes, which serve \revision{as approximations to our scheme} introduced in \eqref{eq:ode_disc}.
The evaluation of the estimator can be another bottleneck in the implementation of the algorithm. For the \revision{higher-order} schemes and \revision{affine right-hand sides}, we noticed that the residual is indeed locally equivalent to the error estimator used in the analysis. \revision{Therefore, one can use} results in \cite[Section 8]{axioms} to justify using the residual itself as estimator, which is accessible \revision{without the need for} $p$th order derivatives of $y_\TT$. {\color{black} This motivates the choice of the residual as an error estimator also for more general non-affine right-hand sides.} \revision{In practice, it might} be advantageous to use an equivalent estimator \revision{based on confidence, namely,} for $T_i=[t_i,t_{i+1}]$,
\begin{align}\label{eq:confidence_estimator}
    \widetilde{\eta}_\TT(T_i)^2 := \frac{\eta_\TT(T_i)^2}{1+\sum_{j=0}^i \eta_\TT(T_j)^2},
\end{align}
or
\begin{align}\label{eq:confidence_estimator_inf}
    \widetilde{\eta}_\TT^\infty(T_i) := \frac{\eta_\TT^\infty(T_i)}{1+\max_{j\leq i} \eta_\TT^\infty(T_j)},
\end{align}
which is more robust. Given a uniform upper bound on the original estimator, one can see that \eqref{eq:confidence_estimator} or \eqref{eq:confidence_estimator_inf} \revision{is locally equivalent to the corresponding original estimator. Thus, the results} in \cite[Section 8]{axioms} justify its use.
\section{Application to neural ODEs}
\revision{We now turn to applications of the adaptive time-stepping framework, beginning with neural ODEs.}

The connection between deep neural networks and ODEs was first made by \cite{Ee}, picked up again by \cite{Haber_2017} and thoroughly investigated by \cite{neuralODE}. \revision{Neural networks built by composing a sequence of hidden states take the form}
\begin{align}
\label{eq:resnet}
    z_{i+1}=z_i+ f_i(z_i, \theta),
\end{align}
\revision{where $i=0,\ldots,N-1$, $z_i \in \R^d$, $\theta$ denotes the weights, and $f_i$ is the composition of an activation function and a linear transformation.} $z_0$ is the input and $z_N$ the output. \revision{The recursion \eqref{eq:resnet} can be seen} as a forward-Euler discretization of some ODE
\begin{align}
\label{eq:neural_ode}
\begin{split}
    \partial_t z (t) &=  f(t,z(t),\theta),\\
    z(t_0) &= z_0,
\end{split}
\end{align}
on some \revision{time interval} $[t_0, \tend]$. \revision{Instead of evaluating the layer structure, we can} parameterise $f$ and define the output of the neural net to be the solution of some \revision{ODE solver} at time $\tend$. \revision{As laid out in \cite{neuralODE}, this} can have several benefits. \revision{One of them is that, during training, there is no need to backpropagate} through the solver, which is usually \revision{very demanding in terms of computation time and memory}. \revision{As discussed below, we can obtain the loss gradient} \revision{with respect to the weights by a single call} to an \revision{ODE solver}, which can be optimized in computational effort by applying the adaptivity results of the previous sections. Many examples use time-independent right-hand sides $f$, which makes them \revision{suitable for our analysis}. \revision{Therefore, in the following, we use} the dynamics
\begin{align}
\label{eq:neural_ode_auto}
\begin{split}
    \partial_t z (t) &=  f(z(t),\theta),\\
    z(t_0) &= z_0.
\end{split}
\end{align}
\subsection{The \revision{adjoint sensitivity method}}
\revision{To apply the adaptive framework during training, we also need an ODE formulation of the loss gradient.}

We consider the neural ODE given by \eqref{eq:neural_ode_auto}, where $z$ is evaluated at $\tend$ as output. A \revision{scalar-valued loss function} $L$ is introduced which depends on $z(\tend)$, i.e.
\begin{align*}
\Loss(z(\tend))=\Loss\Big(z(t_0)+\int_{t_0}^{\tend}f(z(t),\theta)\,\mathrm{d}t\Big).
\end{align*}
In training, we need first derivatives of $L$ with respect to the weights $\theta$. This can be done in a \revision{memory-efficient} way by the \revision{so-called} \textit{adjoint sensitivity method} \cite{neuralODE}. \revision{By formulating the augmented system}
\begin{align}
    \label{eq:ode_aug}
    \begin{split}
        \partial_t z(t) &=  f(z(t),\theta) \\
        \partial_t a(t) &=  - \Big(\nabla_z f(z(t),\theta)\Big)^T a(t), \\
        \partial_t I &=  - \Big(\nabla_{\theta}f(z(t),\theta)\Big)^T a(t),\\
    \end{split}
\end{align}
\revision{and solving it backwards in time with terminal conditions} 
\begin{align*}
    z(\tend) &= z_{\rm end},\\
    a(\tend) &= \frac{\partial L}{\partial z(\tend)},\\
    I(\tend) &= 0,
\end{align*}
where $z_{\rm end}$ is precomputed by solving \eqref{eq:neural_ode_auto}, \revision{we can obtain the loss gradient by evaluating} $I(t_0)$. This is more \revision{memory-efficient}, since general ODE solvers only consider \revision{a fixed number of past values for storage}. Backpropagating through an ODE solver would require much more memory, as the intermediate function evaluations of each layer need to be stored. By defining
\begin{align*}
    y:= \begin{pmatrix}
        z \\
         a\\
        I
    \end{pmatrix},
\end{align*}
and
\begin{align*}
    F(y, \theta):=\begin{pmatrix}
        f(z,\theta) \\
        - \Big(\nabla_z f(z,\theta)\Big)^T a\\
         -\Big(\nabla_{\theta}f(z,\theta)\Big)^T a
    \end{pmatrix},
\end{align*}
\revision{we can write the system \eqref{eq:ode_aug}} as
\begin{align}
\label{eq:ode_aug_simp}
\begin{split}
    \partial_t y(t) &=  F(y(t),\theta),\\
    y(\tend) &= \begin{pmatrix}
        z_{\rm end}\\
        \frac{\partial L}{\partial z(\tend)}\\
        0
    \end{pmatrix},
\end{split}
\end{align}
such that it fits into the discussion of the previous section. \revision{Since $z(t),a(t) \in \R^d$ and $\theta \in \R^w$, we have} $y(t) \in \R^{2d+w}$.
\section{Application to regularized dynamical parametric approximation}\label{sec:regdyn}
\revision{We now extend \eqref{eq:ode} to a Hilbert space setting: for $y_0\in \mathcal{H}$, we seek $y(t)\in \mathcal{H}$ such that}
\begin{align}
\label{eq:ode-hilbert}
\begin{split}
    \partial_t y(t) &= F(y(t)), \qquad \text{for } t \in [t_0, \tend]\\
    y(t_0) &= y_0,
\end{split}
\qquad \text{in } \mathcal{H} ,
\end{align}
for a Lipschitz continuous right-hand side $F:\mathcal{H} \to \mathcal{H}$. 

Parametric dynamical integrators approximate the solution of evolution problems via a nonlinear parameterization
\begin{align*}
    u(t)=\Phi(q(t)).
\end{align*}
Here, $\Phi \, \colon \mathcal Q\rightarrow \mathcal H$ is assumed to be \revision{a smooth parametric map from the parameter space $\mathcal Q \subset \mathbb R^M$, with $M$ parameters, into the space of the evolution problem} \eqref{eq:ode-hilbert}. The path of the parameters $q(t)\in \mathcal Q$ is determined by characterizing its derivative $\dot q (t)$, at any time $t\ge 0$, by the regularized linear least-squares problem
\begin{align*}
    \delta(t)^2 := \norm{\dot{u}(t)-F(u(t))}{}^2+\eps(t)^2\norm{\dot{q}(t)}{}^2 \quad \text{is minimal},
\end{align*}
for given $u(t)$ and $q(t)$ and the regularization parameter $\eps: [t_0,\tend] \to (0,\infty)$. We note that $\dot u$ is fully determined by $\dot q$ through the chain rule $\dot{u} = \Phi'(q)[\dot q]$, where $\Phi'(q)$ denotes the Jacobian of the parametric map $\Phi$ at the position $q$. The regularization of the optimization problem is generally necessary, since we do not want to make any assumption on the condition number of the Jacobian $\Phi'(q)$. In the context of scientific machine learning and computational chemistry, these operators are typically rank-deficient and badly conditioned \cite{regdyn}.

\subsection{An error estimator for the explicit Euler method}
Let $\tend>0$. Introduce a partition $\TT=\{[t_0,t_1],\ldots,[t_{\# \TT-1},\tend]\}$. Compute an approximation $q_\TT \in \SS^1(\TT)$ to $q$ by the explicit Euler method, \revision{i.e.,} the scheme \eqref{eq:ode_disc} for \revision{lowest-order} $p=1$. This means for $t \in T=[t_i,t_{i+1}]$ it holds that
\begin{align*}
    q_{\TT}(t) = q_\TT(t_i) + (t-t_i) \dot{q}_i,
\end{align*}
 where $\dot{q}_i$ solves
\begin{align}\label{eq:regdyn-euler}
    \delta_i^2:=\norm{\Phi'(q_\TT(t_i))\dot{q}_i-F(\Phi(q_\TT(t_i)))}{}^2+\eps_i^2\norm{\dot{q}_i}{}^2 \quad \text{is minimal}.
\end{align}
We note that we use the notation $\|\cdot \|$ \revision{both for the Hilbert space norm in the first summand (typically the $L^2$-norm on some geometry) and for the norm} of the parameter space $\mathcal Q$ in the second summand. {\color{black} One could use the error estimator defined in Section~\ref{sec:low_order} or Section~\ref{sec:low_order_inf} for the parameters $q$, to solve this problem adaptively. As the quantity of interest is an approximation to $y$, we want to motivate a different strategy.} \revision{To this end, we define} another piecewise affine approximation $u_\TT \in \SS^1(\TT)$ to $u$ with $u_\TT(t_0)=y_0$ and
\begin{align}\label{eq:sol_regdyn}
    u_\TT(t) = \Phi(q_\TT(t_i)) + (t-t_i)\frac{\Phi(q_{\TT}(t_{i+1}))-\Phi(q_\TT(t_i))}{|T|}, \quad t\in T.
\end{align} 
\jn{Further, we} define the local error estimator based on $u_\TT$ \jn{by}
\begin{align}\label{eq:esti_regdyn}
    \eta_\TT(T):=|T| \norm{\partial_t u_\TT}{L^2(T)}.
\end{align}
{\color{black}Here $L^2(T)$ denotes the canonical Bochner norm of $L^2(T;\mathcal{H})$. \revision{In the following discussion, the same applies} to $H^1(t_0,\tend)$ norms.}
\revision{One can apply Algorithm~\ref{alg:adaptive} with} the estimator \eqref{eq:esti_regdyn} and solutions obtained via \eqref{eq:sol_regdyn}, \revision{to obtain adaptive approximations to~\eqref{eq:ode-hilbert}} via the regularized  problem. The theory established in this work is not applicable directly to this setting. \revision{Nevertheless, we establish} a related result in the following.
\begin{lemma}\label{lem:rel_regdyn}
    \revision{If $\Phi'$ is Lipschitz continuous with constant $L_2$ and the step-size restriction $|T_i| \delta_i \leq C \eps(t_i)^2$ holds for $T_i=[t_i,t_{i+1}]$, then}
    \begin{align*}
        \norm{y-u_\TT}{H^1(t_0,\tend)}^2\leq C_1 \sum_{T \in \TT} \eta_\TT(T)^2+ C_2\sum_{i=0}^{\# \TT-1} |T_i| \delta_i^2.
    \end{align*}
\end{lemma}
\begin{proof}
    Matching initial conditions and \revision{Gronwall's lemma} give
    \begin{align*}
        \norm{y-u_\TT}{H^1(t_0,\tend)}^2\leq C \norm{\partial_t u_\TT-F(u_\TT)}{L^2(t_0,\tend)}^2.
    \end{align*}
    Looking at $T_i=[t_i,t_{i+1}] \in \TT$ one sees
    \begin{align*}
        \norm{\partial_t u_\TT-F(u_\TT)}{L^2(T_i)}^2 \leq 2\norm{\partial_t u_\TT -F(u_\TT(t_i))}{L^2(T_i)}^2 + 2\norm{F(u_\TT(t_i))-F(u_\TT)}{L^2(T_i)}^2.
    \end{align*}
    \revision{Using the Lipschitz continuity of $F$, we can bound the second term by}
    \begin{align*}
        2L^2|T_i|^2\norm{\partial_t u_\TT}{L^2(T_i)}^2.
    \end{align*}
    For the first term we have by Taylor expansion using the Lipschitz continuity of $\Phi'$ that
    \begin{align*}
        \norm{\partial_t u_\TT -F(u_\TT(t_i))}{L^2(T_i)}^2 \lesssim |T_i|\Big(\norm{\Phi'(q_\TT(t_i))\dot{q}_i-F(u_\TT(t_i))}{}+|T_i| \norm{\dot{q}_i}{}^2\Big)^2.
    \end{align*}
    \revision{The step-size restriction} $|T_i| \delta_i \leq C \eps(t_i)^2$ and the boundedness of $\norm{\Phi'(q_\TT(t_i))\dot{q}_i-F(u_\TT(t_i))}{}$ \revision{by $\delta_i$ give}
    \begin{align*}
        \norm{\partial_t u_\TT -F(u_\TT(t_i))}{L^2(T_i)}^2 \lesssim |T_i| \delta_i^2.
    \end{align*}
    \revision{Summing over all $T_i$} concludes the proof.
\end{proof}
Lemma~\ref{lem:rel_regdyn} \revision{establishes that} the proposed estimator is \revision{reliable up to the approximation error} of the regularized parametric problem. We denote by $\TT_0,\TT_1,\ldots$ a sequence of meshes and solutions obtained by Algorithm~\ref{alg:adaptive} to state the next result. Superscripts are used to signal that the discrete defects $\delta_i^k$ are associated with mesh $\TT_k$. 
\begin{lemma}\label{lem:ada_regdyn} \revision{If $\Phi'$ is Lipschitz continuous with constant $L_2$ and the step-size restriction} $|T_i| \delta_i \leq C \eps(t_i)^2$ holds for $T_i=[t_i,t_{i+1}] \in \TT_k$, for all meshes $\TT_k$, then
    there exist constants $C_{\rm stab},C_{\rm red}>0$ and $q \in (0,1)$ such that for all $k \in \mathbb{N}_0$ it holds that
     \begin{align}\label{eq:stab_regdyn}
     \begin{split}
          \Bigg| \Big(\sum_{T \in \TT_k \cap \TT_{k+1}}\eta_{\TT_k}(T)^2 \Big)^{1/2} &-\Big(\sum_{T \in \TT_k \cap \TT_{k+1}} \eta_{\TT_{k+1}}(T)^2\Big)^{1/2} \Bigg|
        \\
 &\leq  C_{\rm stab} \norm{h_{\TT_k}}{L^\infty(t_0,\tend)}\norm{u_{\TT_k}-u_{\TT_{k+1}}}{H^1(t_0,\tend)}
     \end{split}
    \end{align}
    \begin{align}\label{eq:red_regdyn}
    \begin{split}
        \sum_{T \in \TT_{k+1} \setminus \TT_k} \eta_{\TT_{k+1}}(T)^2 &\leq  q \sum_{T \in  \TT_k\setminus \TT_{k+1}} \eta_{\TT_k}(T)^2 \\
 &+ C_{\rm red}\norm{h_{\TT_k}}{L^\infty(t_0,\tend)}^2\norm{u_{\TT_k}-u_{\TT_{k+1}}}{H^1(t_0,\tend)}^2.
    \end{split}
    \end{align}
\revision{Furthermore, for some $0<q_{\rm red}<1$, we have}
\begin{align}\label{eq:est_red_regdyn}
\eta_{\TT_{k+1}}^2 \leq q_{\rm red} \eta_{\TT_{k}}^2 + C \norm{h_{\TT_k}}{L^\infty(t_0,\tend)}^2\norm{u_{\TT_k}-u_{\TT_{k+1}}}{H^1(t_0,\tend)}^2.
\end{align}
\revision{For sufficiently small $\norm{h_{\TT_0}}{L^\infty(t_0,\tend)}$, the estimator $\eta_{\TT_k}$ contracts up to a perturbation term:}
\begin{align}\label{eq:lin_conv_regdyn1}
    \eta_{\TT_{k+1}}^2 \leq \kappa \eta_{\TT_k}^2 + C_3 \Bigg( \sum_{i=0}^{\# \TT_k-1} |T_i^{k}| (\delta^k_i)^2 + \sum_{i=0}^{\# \TT_{k+1}-1} |T_i^{k+1}|(\delta^{k+1}_i)^2\Bigg),
\end{align}
with 
\begin{align*}
    \kappa &:=\frac{q_{\rm red}+ CC_1\norm{h_{\TT_0}}{L^\infty(t_0,\tend)}^2}{1-CC_1\norm{h_{\TT_0}}{L^\infty(t_0,\tend)}^2}<1,\\
    C_3 &:=\frac{CC_2\norm{h_{\TT_0}}{L^\infty(t_0,\tend)}^2}{1-CC_1\norm{h_{\TT_0}}{L^\infty(t_0,\tend)}^2}.
\end{align*}
\revision{Writing $\delta^k:=\max_i \delta^k_i$, one obtains} for $N \in \mathbb{N}$ that
\begin{align}\label{eq:lin_conv_regdyn2}
    \eta_{\TT_{k+N}}^2 \leq \kappa^N \eta_{\TT_k}^2 + C_3 (\tend-t_0)\Big( \sum_{j=0}^{N-1} (\delta^{k+N-j})^2 \kappa^{j}+ \sum_{j=1}^{N} (\delta^{k+N-j})^2 \kappa^{j-1}\Big).
\end{align}
\end{lemma}
\begin{proof}
    \eqref{eq:stab_regdyn}, \eqref{eq:red_regdyn} and \eqref{eq:est_red_regdyn} follow \jn{\revision{by arguments analogous to those presented} in Section~\ref{sec:low_order}.} For \eqref{eq:lin_conv_regdyn1} we use \eqref{eq:est_red_regdyn} combined with Lemma~\ref{lem:rel_regdyn} to obtain
    \begin{align*}
        \eta_{\TT_{k+1}}^2 &\leq  (q_{\rm red}+CC_1\norm{h_{\TT_k}}{L^\infty(t_0,\tend)}^2) \eta_{\TT_k}^2+  CC_1 \norm{h_{\TT_{k}}}{L^\infty(t_0,\tend)}^2 \eta_{\TT_{k+1}}^2\\
 &+CC_2 \norm{h_{\TT_{k}}}{L^\infty(t_0,\tend)}^2 \Bigg( \sum_{i=0}^{\# \TT_k-1} |T_i^{k}| (\delta^k_i)^2 + \sum_{i=0}^{\# \TT_{k+1}-1} |T_i^{k+1}|(\delta^{k+1}_i)^2\Bigg).
    \end{align*}
    \revision{Choosing a sufficiently fine initial mesh} such that $CC_1 \norm{h_{\TT_{k}}}{L^\infty(t_0,\tend)}^2 \eta_{\TT_{k+1}}^2$ can be shifted to the other side and such that $\kappa<1$ yields the result. \eqref{eq:lin_conv_regdyn2} follows by consecutive application of \eqref{eq:lin_conv_regdyn1}.
\end{proof}
The interpretation of Lemma~\ref{lem:ada_regdyn} is the following. Provided that the defects admit a uniform bound independent of the underlying mesh, the adaptive algorithm reduces the error estimator up to this defect bound. This is exactly the desired behavior in the present application.
{\color{black}\subsection{Generalization to \revision{higher-order} schemes}\label{sec:regdyn_higher}
The generalization of the error estimator to the \revision{higher-order} schemes introduced in Sections~\ref{sec:high_order} \revision{and~\ref{sec:high_order_inf} is not immediately clear.} The preceding section provides some motivation to take an error estimator based on the residual of $\Phi(q_\TT)$. We will investigate the performance of this choice in the numerical experiments.}
\section{Numerical Experiments}
\revision{We now illustrate the adaptive schemes for the model problems covered by the analysis and for the applications discussed above.}

{\color{black}In the following experiments, the estimator is chosen according to the problem and the discretization order. Whenever the assumptions of the analysis are satisfied, we use the corresponding estimator introduced above: \eqref{eq:estimator} and \eqref{eq:esti_inf} for the first-order scheme, and \eqref{eq:lin_esti} and \eqref{eq:lin_esti_inf} for higher-order schemes with affine right-hand sides. The \(L^2\)-based indicators are combined with Dörfler marking, while the \(L^\infty\)-based indicators are combined with maximum marking.
For nonlinear problems treated with higher-order schemes, which lie outside the scope of the optimality analysis, we use the residual \(R=\partial_t y_{\mathcal T}-F(y_{\mathcal T})\) to define the local indicators
\[
\eta_{\mathcal T}(T):=\|R\|_{L^2(T)}
\quad \text{\jn{and}} \quad
\eta_{\mathcal T}^{\infty}(T):=\|R\|_{L^\infty(T)}.
\]
This choice is motivated by Lemma~\ref{lem:rel} and the discussion in Section~\ref{sec:practical}. The corresponding experiments therefore investigate the applicability of this approach beyond the proven setting. For regularized dynamical parametric approximation, we use \eqref{eq:esti_regdyn} for the first-order scheme and \jn{higher-order derivatives of \(\Phi(q_{\mathcal T})\) for higher-order schemes, analogous to \eqref{eq:esti_regdyn} and \eqref{eq:lin_esti}.} Where indicated, the chosen indicators are further normalized according to \eqref{eq:confidence_estimator} or \eqref{eq:confidence_estimator_inf}.}
\subsection{Heat equation}
\begin{figure}[htbp]
    \centering

    \begin{subfigure}[t]{0.32\textwidth}
        \centering
        \includegraphics[width=\linewidth,trim=22bp 0bp 37bp 37bp,clip]{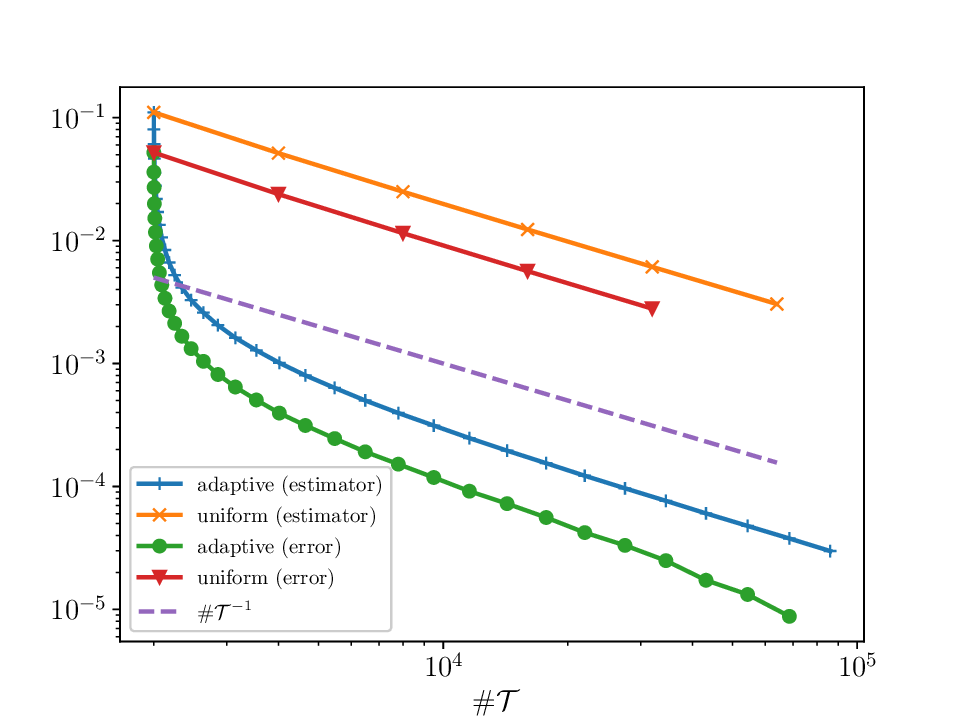}
        \caption{$\theta=0.5$, $\eta_\TT$}
    \end{subfigure}\hspace{0.005\textwidth}%
    \begin{subfigure}[t]{0.32\textwidth}
        \centering
        \includegraphics[width=\linewidth,trim=22bp 0bp 37bp 37bp,clip]{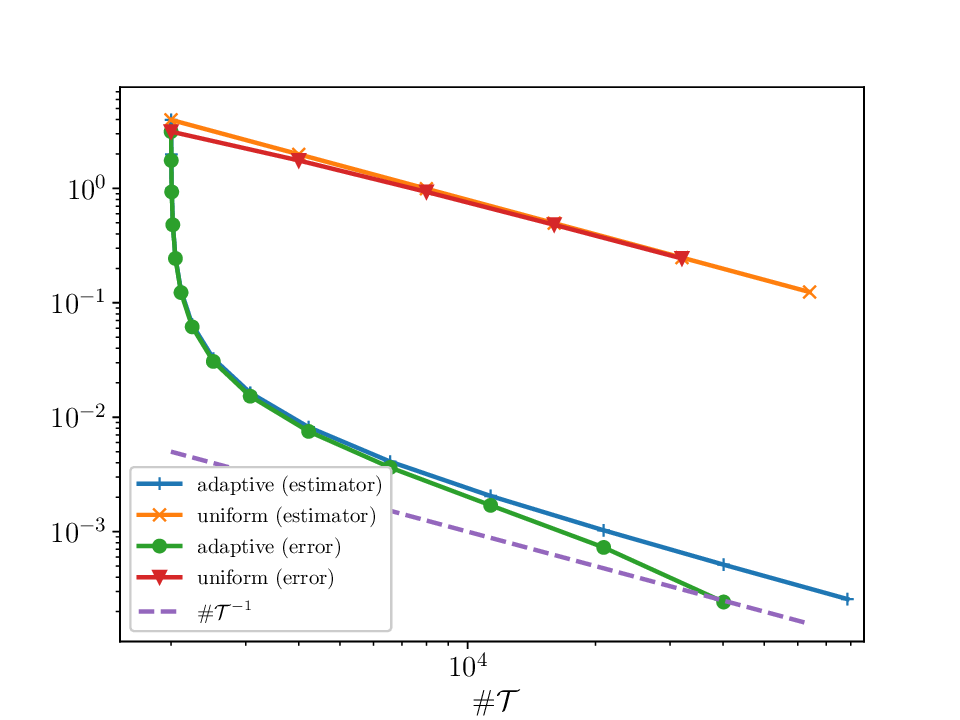}
        \caption{$\theta=0.5$, $\eta_\TT^\infty$}
    \end{subfigure}\hspace{0.005\textwidth}%
    \begin{subfigure}[t]{0.32\textwidth}
        \centering
        \includegraphics[width=\linewidth,trim=22bp 0bp 37bp 37bp,clip]{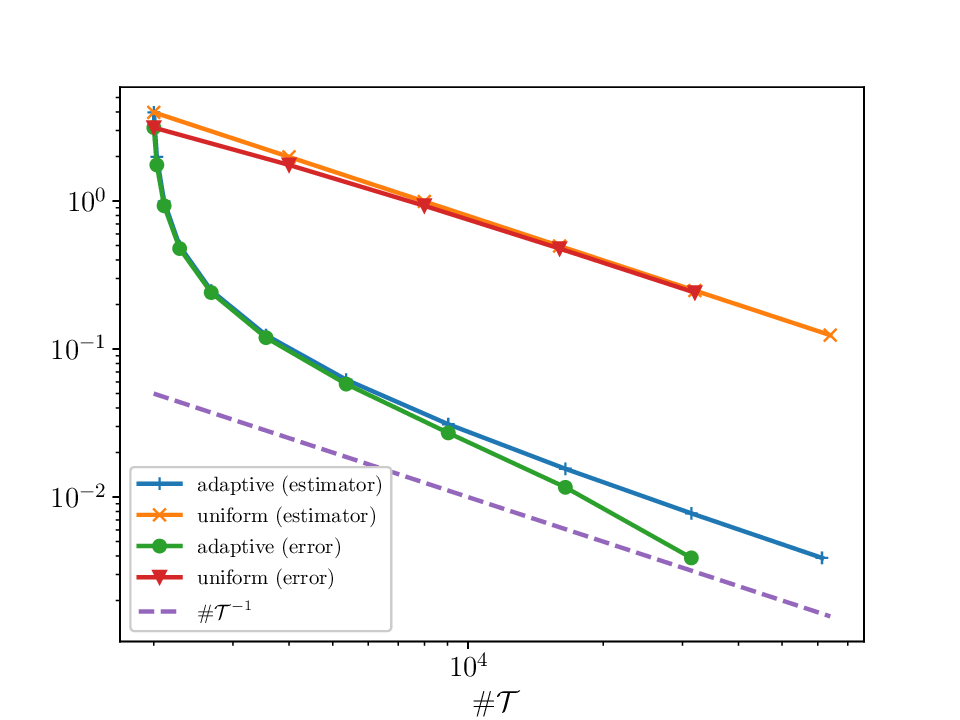}
        \caption{$\theta=0.99$, $\eta_\TT^\infty$}
    \end{subfigure}

    \par\smallskip

    \begin{subfigure}[t]{0.32\textwidth}
        \centering
        \includegraphics[width=\linewidth,trim=22bp 0bp 37bp 37bp,clip]{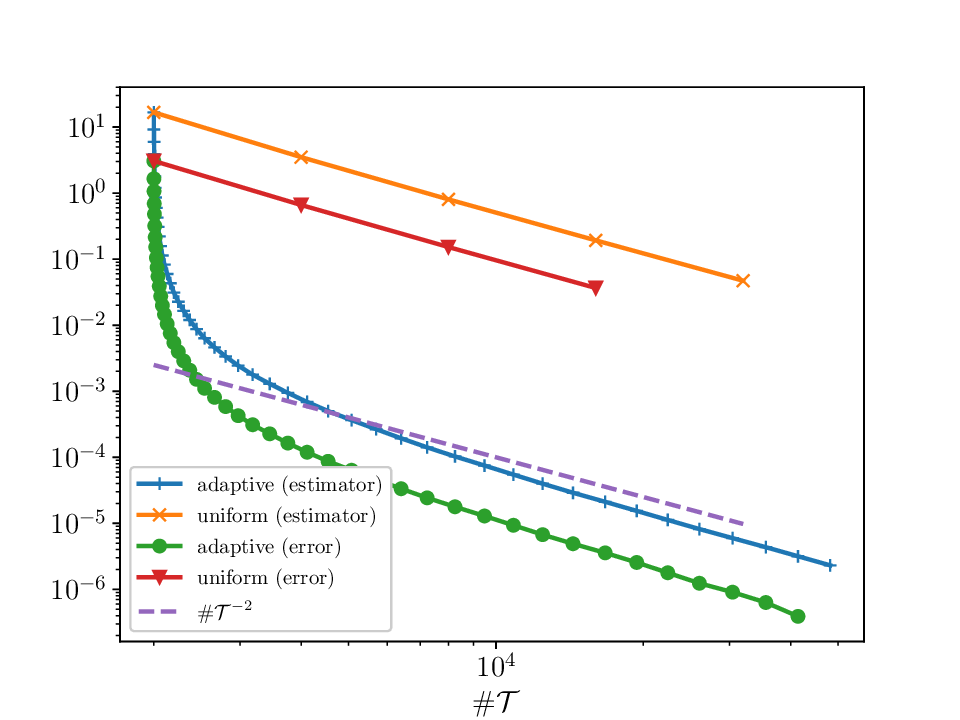}
        \caption{$\theta=0.5$, $\eta_\TT$}
    \end{subfigure}\hspace{0.005\textwidth}%
    \begin{subfigure}[t]{0.32\textwidth}
        \centering
        \includegraphics[width=\linewidth,trim=22bp 0bp 37bp 37bp,clip]{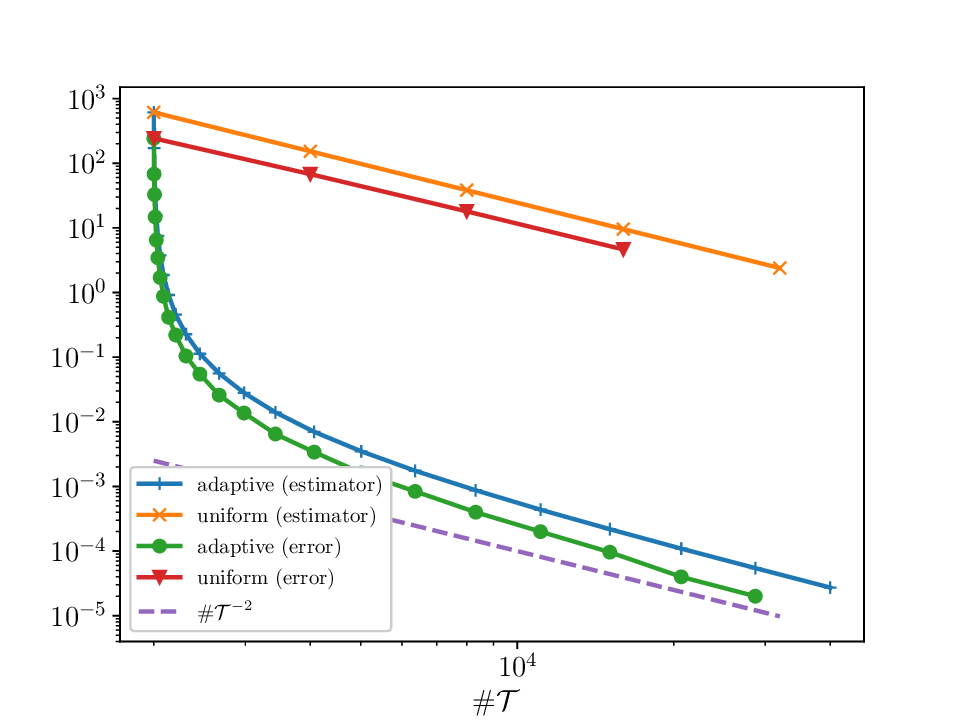}
        \caption{$\theta=0.5$, $\eta_\TT^\infty$}
    \end{subfigure}\hspace{0.005\textwidth}%
    \begin{subfigure}[t]{0.32\textwidth}
        \centering
        \includegraphics[width=\linewidth,trim=22bp 0bp 37bp 37bp,clip]{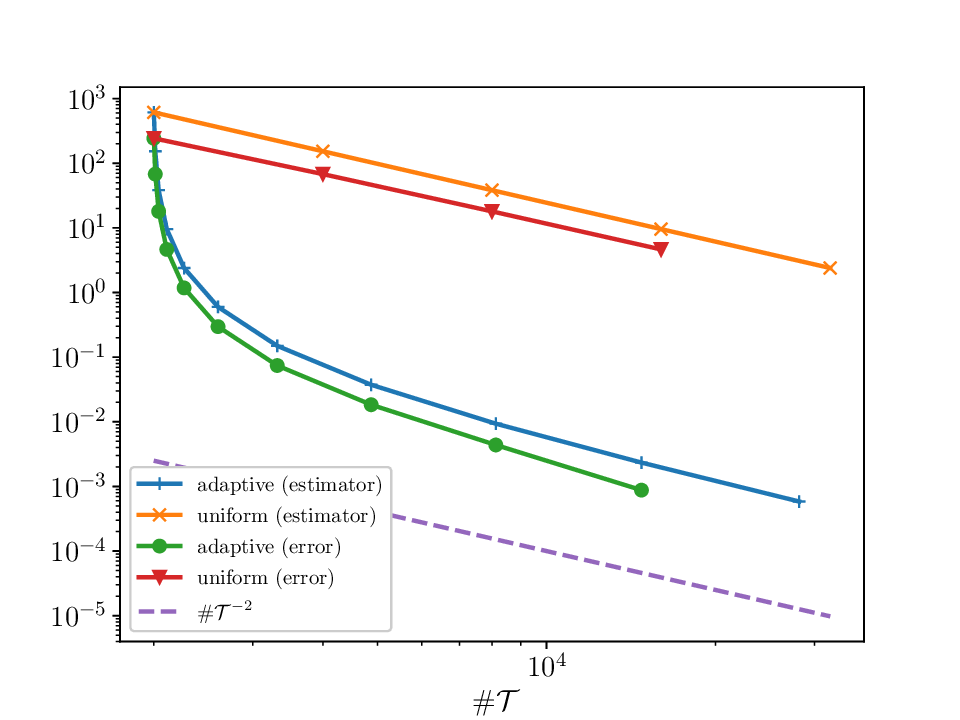}
        \caption{$\theta=0.99$, $\eta_\TT^\infty$}
    \end{subfigure}

    \caption{\color{black}
        Heat equation: estimators and corresponding errors versus the
        number of time intervals $\#\mathcal{T}$ for adaptive and uniform
        refinement. The top row (a)--(c) shows results for $p=1$, and
        the bottom row (d)--(f) shows results for $p=2$.
        Adaptive refinement uses Dörfler marking in (a) and (d),
        and maximum marking in (b)--(c) and (e)--(f).
        The dashed lines indicate the reference rates
        $(\#\mathcal{T})^{-p}$.
    }
    \label{fig:heat}
\end{figure}

\begin{figure}[h]
    \centering
    \begin{subfigure}[t]{0.49\textwidth}
    \centering
    \includegraphics[width=\linewidth,trim=22bp 0bp 39bp 38bp,clip]{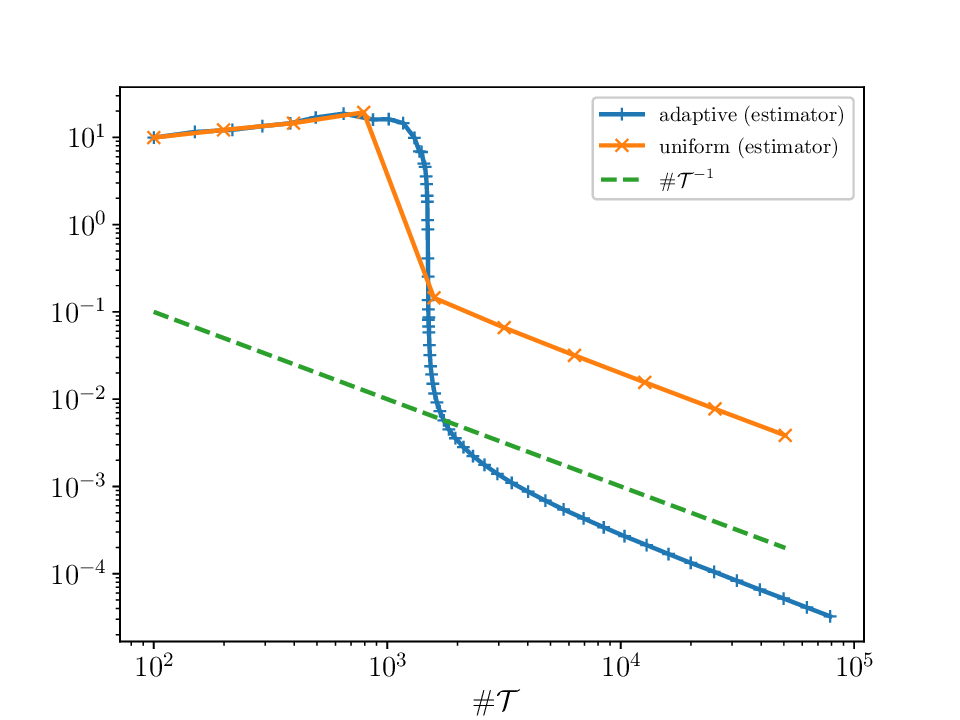}
    \caption{$p=1$, $\widetilde\eta_\TT$}
\end{subfigure}\hspace{0.005\textwidth}%
\begin{subfigure}[t]{0.49\textwidth}
    \centering
    \includegraphics[width=\linewidth,trim=22bp 0bp 39bp 38bp,clip]{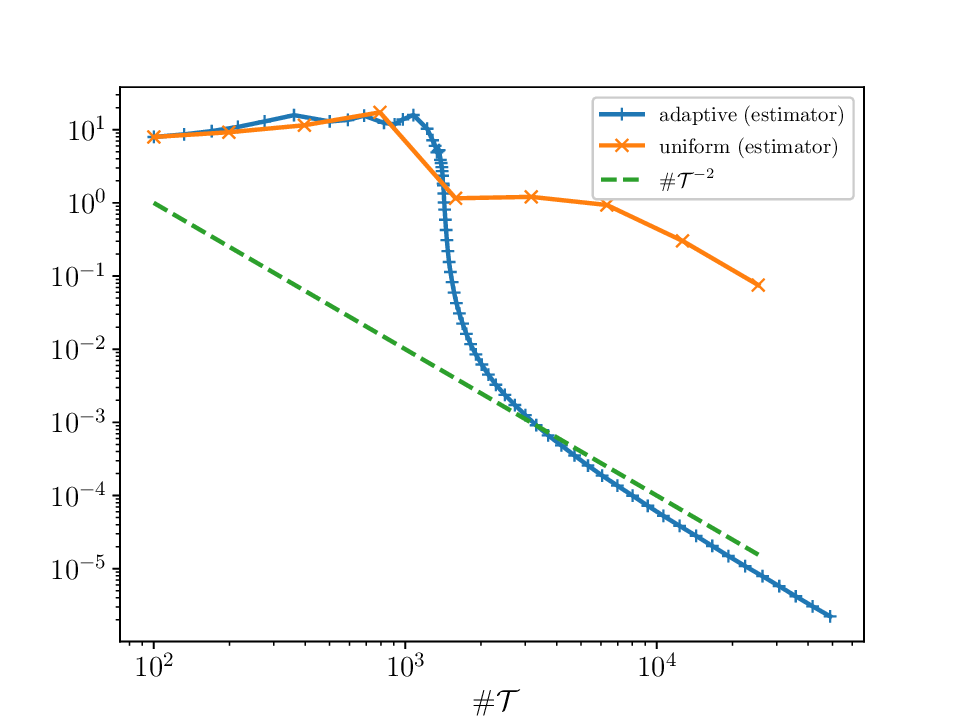}
    \caption{$p=2$, $\widetilde\eta_\TT$}
\end{subfigure}

\par\smallskip

 \begin{subfigure}[t]{0.49\textwidth}
    \centering
    \includegraphics[width=\linewidth,trim=22bp 0bp 39bp 38bp,clip]{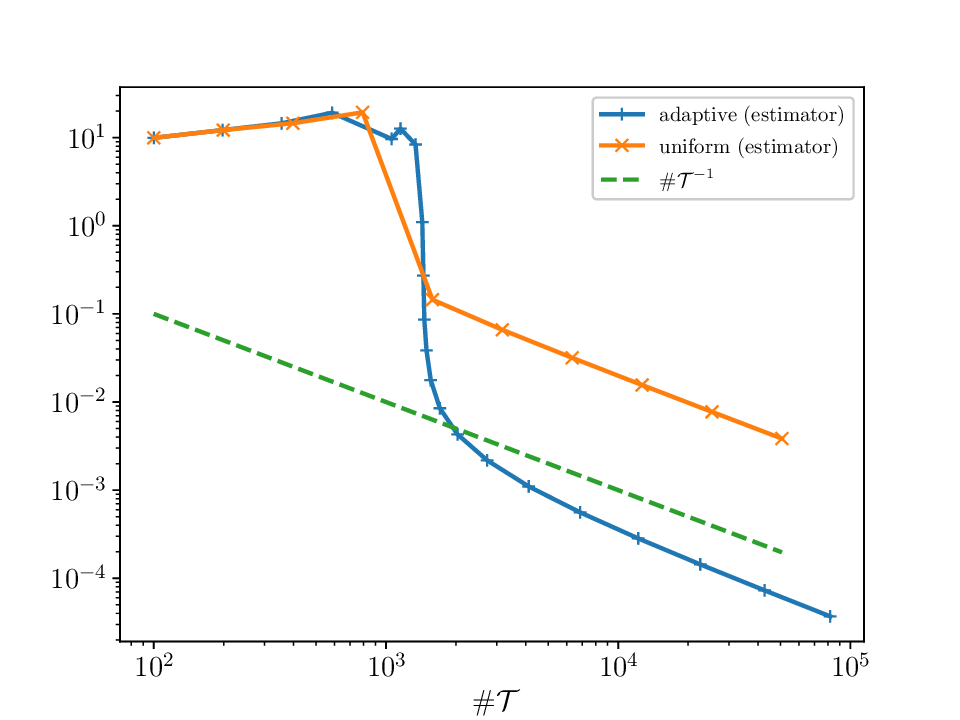}
    \caption{$p=1$, $\widetilde\eta_\TT^\infty$}
\end{subfigure}\hspace{0.005\textwidth}%
\begin{subfigure}[t]{0.49\textwidth}
    \centering
    \includegraphics[width=\linewidth,trim=22bp 0bp 39bp 38bp,clip]{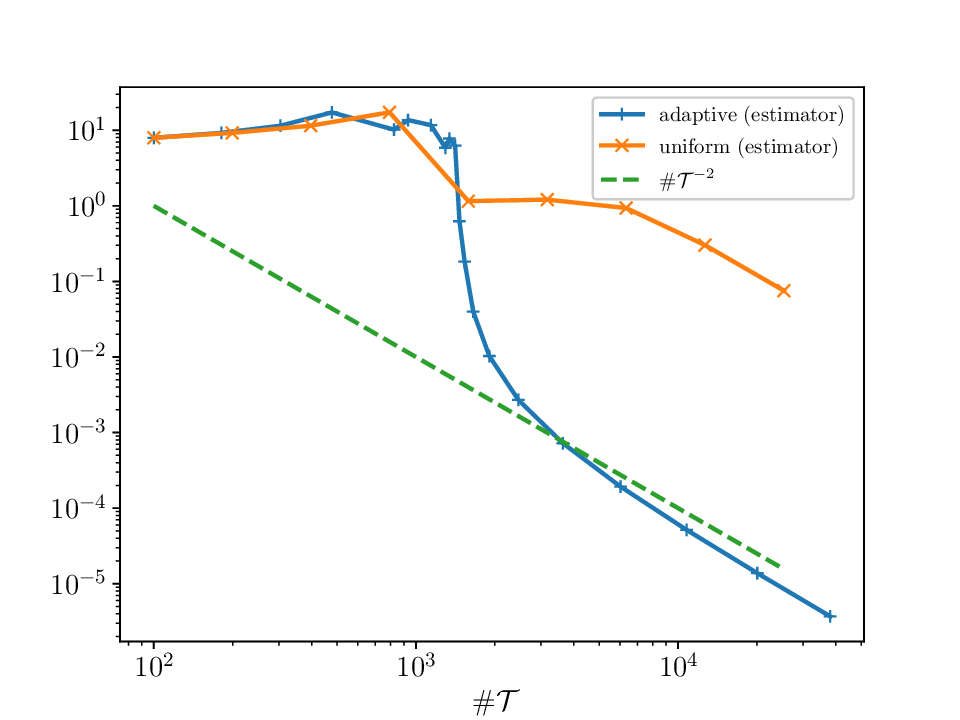}
    \caption{$p=2$, $\widetilde\eta_\TT^\infty$}
\end{subfigure}
\caption{\color{black} Heat equation: confidence-normalized indicators versus the number of time intervals $\#\mathcal{T}$ for adaptive and uniform refinement, with orders $p=1,2$. Adaptive refinement uses Dörfler marking in (a)--(b) and maximum marking in (c)--(d). The dashed lines indicate the reference rates $(\#\mathcal{T})^{-p}$.}
 \label{fig:confidence}
\end{figure}

\revision{Our first benchmark example is} the heat equation,
\begin{align*}
    \partial_t u - \Delta u &=  f, \quad \text{on } [0,1] \times \Omega,\\
    u(0) &= u_0.
\end{align*}
with $\Omega=(0,1)^2$, and homogeneous Dirichlet boundary conditions on $\partial \Omega$. We semidiscretize this problem \revision{in space by} a conforming P1 finite element space on a uniform triangulation of $\Omega$, with \revision{mesh width} $h=1/10$. This yields a linear ODE, which fits into the discussion above. \revision{We choose the right-hand side $f=0$, take $u_0$ as the $L^2(\Omega)$-projection of the constant function $1$} onto the chosen finite element space and $\tend=1$. We apply our scheme \eqref{eq:ode_disc} with $p=1$ and $p=2$, $\theta \in \{0.5,0.99\}$ and use both the Dörfler marking and maximum strategy for refinement with the corresponding error estimators. {\color{black}Figure~\ref{fig:heat} compares adaptive and uniform temporal refinement for \(p=1\) and \(p=2\). The adaptive results exhibit decay consistent with the reference rates \((\#\mathcal T)^{-p}\) and attain smaller estimator and error values at comparable numbers of time intervals. Figure~\ref{fig:confidence} presents the corresponding experiments with the confidence-normalized indicators. Starting from coarse temporal meshes, the adaptive indicators initially exhibit a plateau, followed by decay consistent with the respective reference rates. The actual errors are omitted, as they cannot be meaningfully displayed during the preasymptotic phase.}

\subsection{Van-der-Pol equation: influence of $\mu$}
\begin{figure}[h!]
    \centering
    \begin{subfigure}[t]{0.32\textwidth}
        \centering
        \includegraphics[width=\linewidth,trim=22bp 0bp 42bp 35bp,clip]{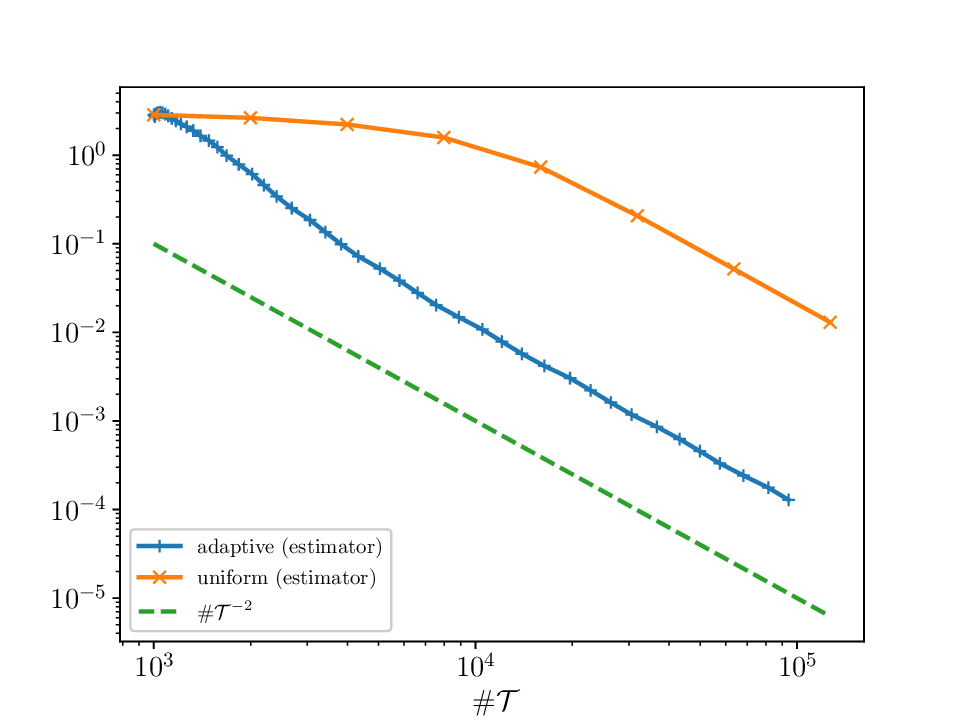}
        \caption{$p=2$, $\mu=100$, $\widetilde{\eta}_\TT$}
    \end{subfigure}\hspace{0.005\textwidth}%
    \begin{subfigure}[t]{0.32\textwidth}
        \centering
        \includegraphics[width=\linewidth,trim=22bp 0bp 42bp 35bp,clip]{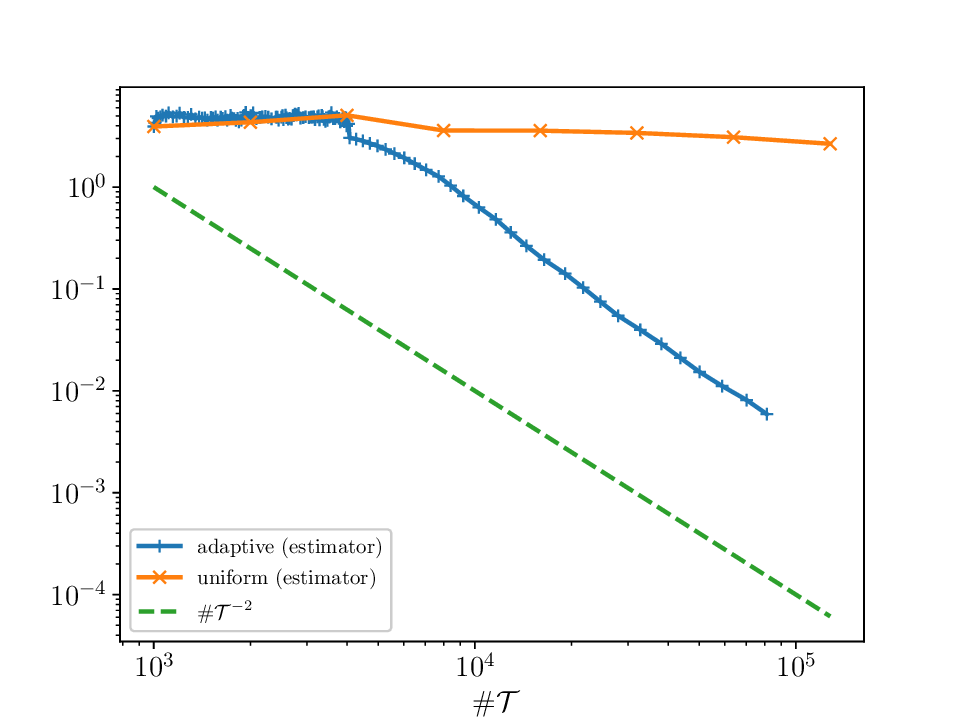}
        \caption{$p=2$, $\mu=1000$, $\widetilde{\eta}_\TT$}
    \end{subfigure}\hspace{0.005\textwidth}%
    \begin{subfigure}[t]{0.32\textwidth}
        \centering
        \includegraphics[width=\linewidth,trim=22bp 0bp 42bp 35bp,clip]{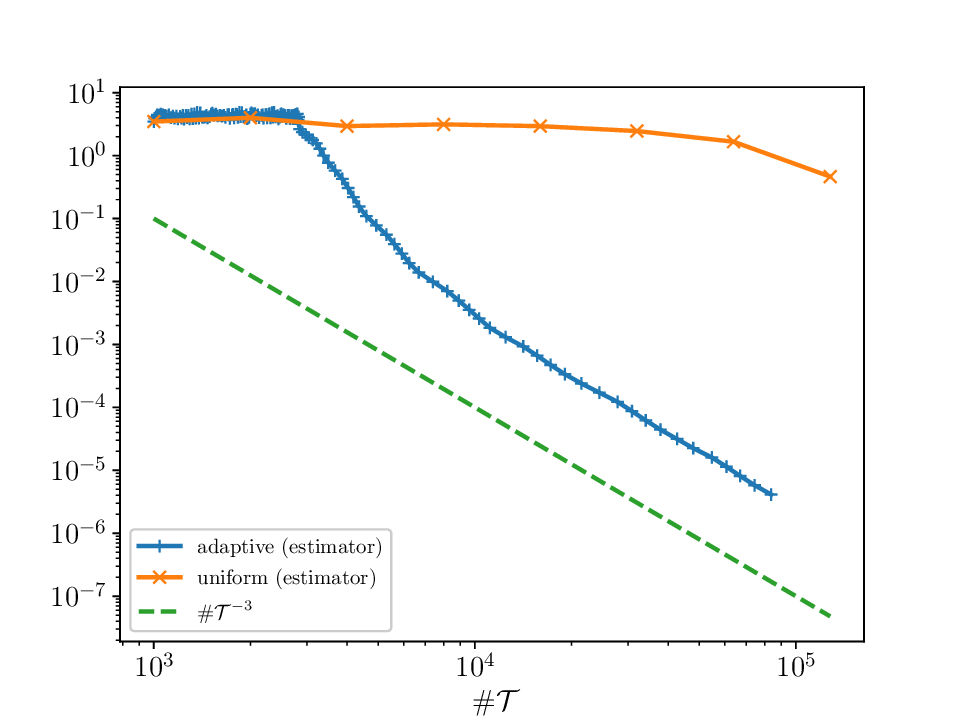}
        \caption{$p=3$, $\mu=1000$, $\widetilde{\eta}_\TT$}
    \end{subfigure}

    \par\smallskip

    \begin{subfigure}[t]{0.32\textwidth}
        \centering
        \includegraphics[width=\linewidth,trim=22bp 0bp 42bp 35bp,clip]{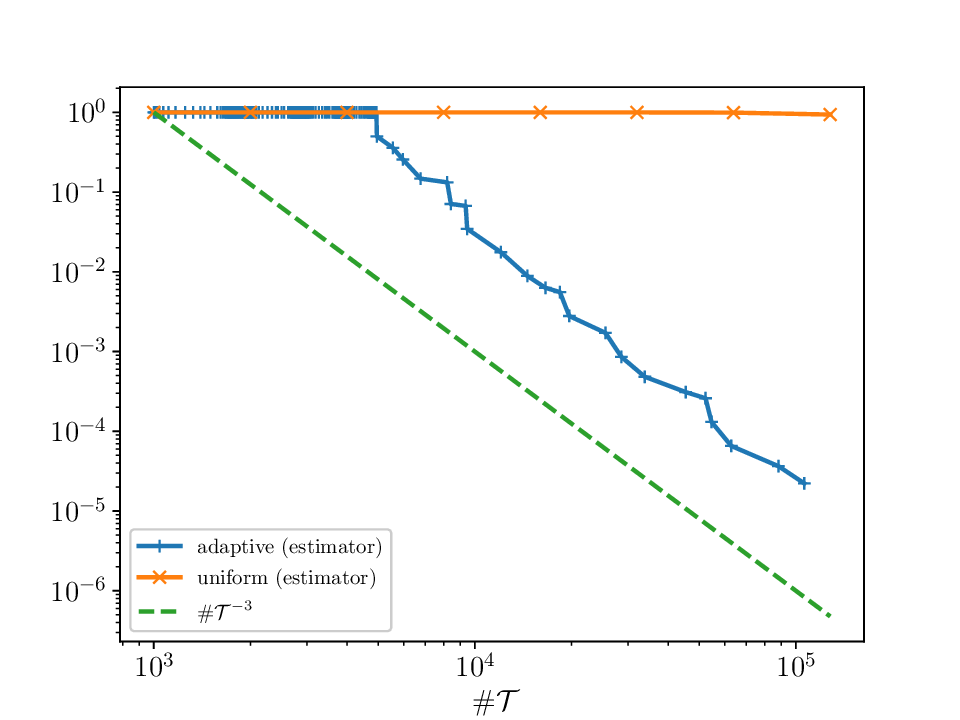}
        \caption{$p=3$, $\mu=1000$, $\widetilde{\eta}_\TT^\infty$}
    \end{subfigure}\hspace{0.005\textwidth}%
    \begin{subfigure}[t]{0.32\textwidth}
        \centering
        \includegraphics[width=\linewidth,trim=22bp 0bp 42bp 35bp,clip]{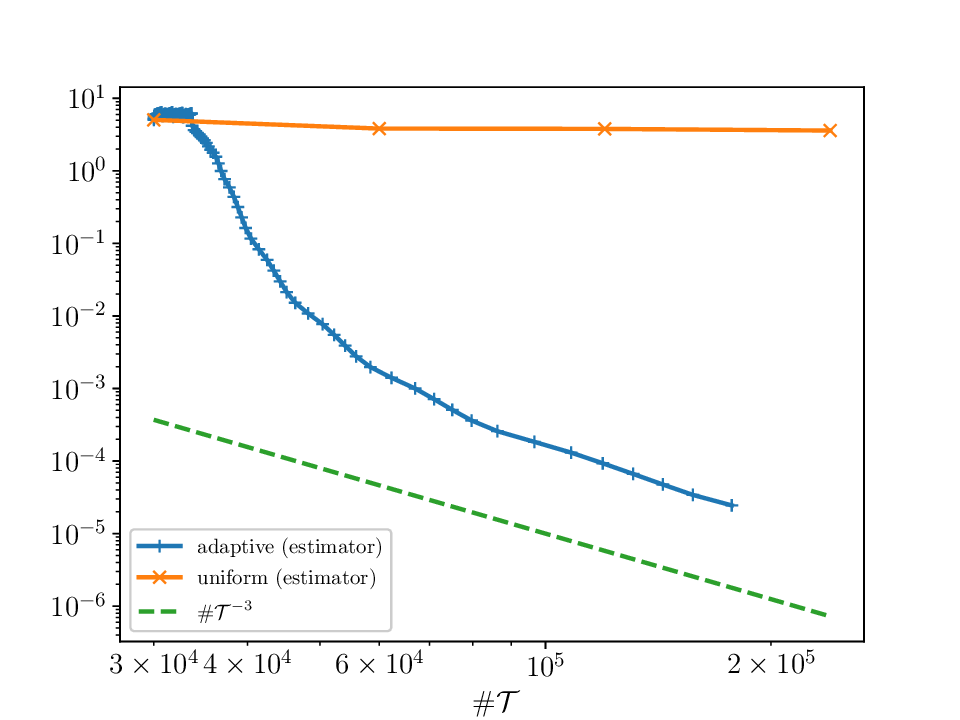}
        \caption{$p=3$, $\mu=10000$, $\widetilde{\eta}_\TT$}
    \end{subfigure}\hspace{0.005\textwidth}%
    \begin{subfigure}[t]{0.32\textwidth}
        \centering
        \includegraphics[width=\linewidth,trim=22bp 0bp 42bp 35bp,clip]{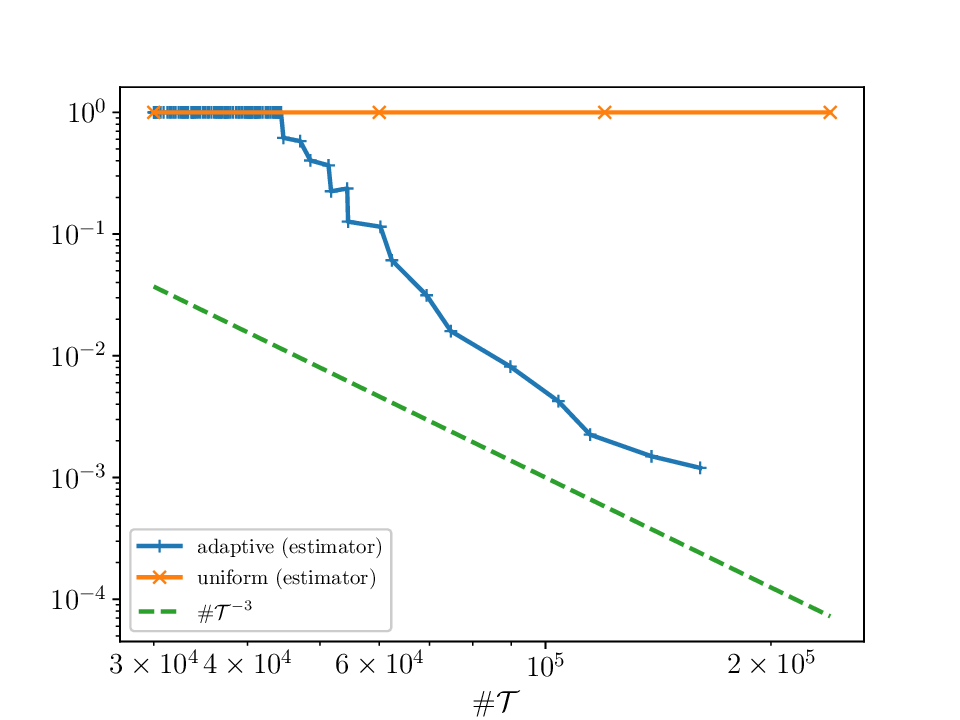}
        \caption{$p=3$, $\mu=10000$, $\widetilde{\eta}_\TT^\infty$}
    \end{subfigure}
    \caption{\color{black} Van der Pol problem: confidence-normalized residual indicators versus the number of time intervals $\#\mathcal{T}$ for adaptive and uniform refinement, with orders $p=2,3$, marking parameter $\theta=0.5$ and varying stiffness parameter $\mu$. Adaptive refinement uses Dörfler marking in (a)--(c) and (e), and maximum marking in (d) and (f). The dashed lines indicate the reference rates $(\#\mathcal{T})^{-p}$.}
    \label{fig:vdp_conf}
\end{figure}

The Van-der-Pol equation is a standard toy problem to test robust \revision{time-stepping method}s. It is stated as
\begin{align*}
     \partial_t x =y, \quad
    \partial_t y = \mu(y(1-x^2)-x),
\end{align*}
where $\mu$ determines the stiffness of the system. We choose $(x(0),y(0))=(1,1)$ 
and test the robustness of the adaptive algorithm equipped with the \revision{confidence-based} estimator with respect to $\mu$. We solve the equation on the interval $[0,3]$. {\color{black}We investigate the influence of \(\mu\) using the second-order scheme with \(\mu\in\{100,1000\}\) and the third-order scheme with \(\mu\in\{1000,10000\}\). Adaptivity is driven by the confidence-normalized residual indicators introduced in Section~\ref{sec:practical}. Figure~\ref{fig:vdp_conf} compares adaptive and uniform refinement, with both Dörfler and maximum marking considered for the third-order scheme.
Over the displayed range, adaptive refinement produces substantially smaller normalized residual indicators than uniform refinement at comparable numbers of time intervals. For larger values of \(\mu\), clear estimator decay begins \revision{at larger numbers of time steps}. The dashed lines indicate the reference rates \((\#\mathcal T)^{-p}\).}
Note that our schemes are explicit and still produce reasonable results. \jn{For large values of $\mu$, the system becomes stiff and explicit methods require a very fine initial grid. \revision{This manifests itself in a growth} of the preasymptotic phase, which can become computationally infeasible for many practical problems. In such cases, one should rely on \revision{implicit adaptive methods to compute accurate solutions with a reasonable number of time steps}.} We refer to the seminal book \cite{hairer_stiff} or \cite{implicitadaptive} for a corresponding algorithm.
\subsection{Transport problem with varying coefficients} 
\begin{figure}[h!]
\centering
\includegraphics[width=1\linewidth]{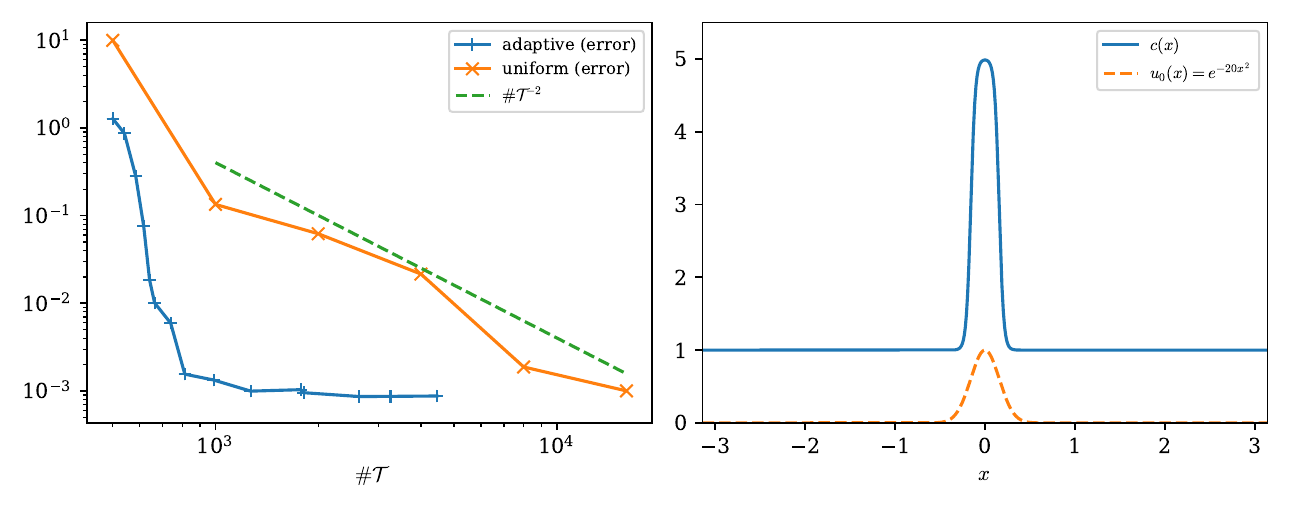}
\caption{On the left: convergence plot for the uniform explicit midpoint rule and the adaptive $p=2$ rule for the transport problem with spatially varying, localized coefficients \eqref{eq:transport}. The plot on the right-hand side shows the coefficient function $c(x)$ and the initial condition.}
\label{fig:err-coeff}
\end{figure}

We consider the one-dimensional transport equation
\begin{equation}\label{eq:transport}
\partial_t y(x,t) = c(x) \partial_x y(x,t), \qquad \text{for } x\in (-\pi,\pi).  
\end{equation}
Along the boundary, we enforce periodic boundary conditions. \jn{We note that this model problem violates the assumed Lipschitz continuity that is required for the optimality of the proposed adaptive time-stepping scheme. The experiment is therefore designed to test the proposed methodology outside of the theory that we provide within this manuscript. Rigorously showing an optimality result for the following adaptive procedure is therefore left for future work.}
We use the regularized dynamical approximation of Section~\ref{sec:regdyn}, which constructs the approximation $\Phi(q(t))$. To observe some challenging dynamics in the solution, we use the \revision{initial condition and coefficient} \begin{align}\label{eq:init-coeff}
\begin{aligned}
y_0(x)&=\exp(-20x^2), \ \ \ \ \ \ \\ c(x)&=1+2\left(\tanh\left(\dfrac{x+0.05\pi}{0.05}\right)-\tanh\left(\dfrac{x-0.05\pi}{0.05}\right)\right),
\end{aligned}
\end{align}
which are both smooth but spike relatively sharply around the origin (see Figure~\ref{fig:err-coeff}). The dynamics of the evolution problem are therefore \revision{relatively simple most of the time}, but challenging when the support of the solution at a time intersects with the region around the origin in which the coefficients differ sharply from the rest of the domain.

As the initial parametrization, we fit a neural network architecture that enforces periodic boundary conditions exactly, uses four fully connected hidden layers with the $\tanh$ activation function and $10$ neurons each, \revision{with $461$ parameters in total}. A more detailed \revision{description of the architecture} and the fitting process can be found in \cite[Section~8.2]{LN25}.

The parameters are then evolved over time by different regularized parametric approaches. 
For all experiments, we set the regularization parameter to $\varepsilon =10^{-3}$.  First, we consider the regularized parametric explicit midpoint method with uniform step sizes. \revision{Second, we consider} the parametric variant of the second order Runge--Kutta scheme of Figure~\ref{fig:rk} with adaptive refinements. We briefly sketch this integrator to define the error indicator used. First, we define
\begin{align*}
    Q_{i,1} &:= q_\TT(t_i), &
    Q_{i,2} &:= q_\TT(t_i)+\tfrac{h_i}{3}\,\dot{q}_{i,1}, &
    Q_{i,3} &:= q_\TT(t_i)+h_i\,\dot{q}_{i,1},
\end{align*}
where $\dot{q}_{i,s}$ for $s=1,2,3$ solves
\begin{align*}
\delta_{i,s}^2:=\norm{\Phi'(Q_{i,s})\,\dot{q}_{i,s}-F(\Phi(Q_{i,s}))}{}^2+\eps^2\norm{\dot{q}_{i,s}}{}^2 \quad \text{is minimal}.
\end{align*}
Then, for $t\in T_i=[t_i,t_{i+1}]$,
\begin{align*}
    q_\TT(t) = q_\TT(t_i)
    + (t-t_i)\bigl(\tfrac32\,\dot{q}_{i,2}-\tfrac12\,\dot{q}_{i,3}\bigr)
    + \frac{(t-t_i)^2}{2}\,\ddot{q}_i,
    \qquad
    \ddot{q}_i := \frac{3}{2h_i}\bigl(\dot{q}_{i,3}-\dot{q}_{i,2}\bigr).
\end{align*}
In order to approximate the error indicator \eqref{eq:lin_esti_inf}, we approximate the second temporal derivative by neglecting curvature terms of $\Phi$, which gives
\begin{align*}
\partial_t^2 \Phi (q_{\TT}(t))
 &\approx
\Phi'(q_{\mathcal T}(t_i))[\ddot{q}_i].
\end{align*}
We then use \revision{the following approximate error indicator}
\begin{align*}
    \eta_\TT^\infty(T_i) = h_i^2\,\norm{F\bigl(\Phi'(q_{\mathcal T}(t_i))\,[\ddot{q}_i]\bigr)}{}.
\end{align*}

We use the maximum marking strategy of Algorithm~\ref{alg:adaptivemax} with the marking parameter $\theta=0.9$. \revision{The initial temporal mesh for adaptive refinement decomposes} the time interval $[0,10]$ into $N=500$ time steps. 

\begin{figure}[h!]
\centering
\includegraphics[width=1\linewidth]{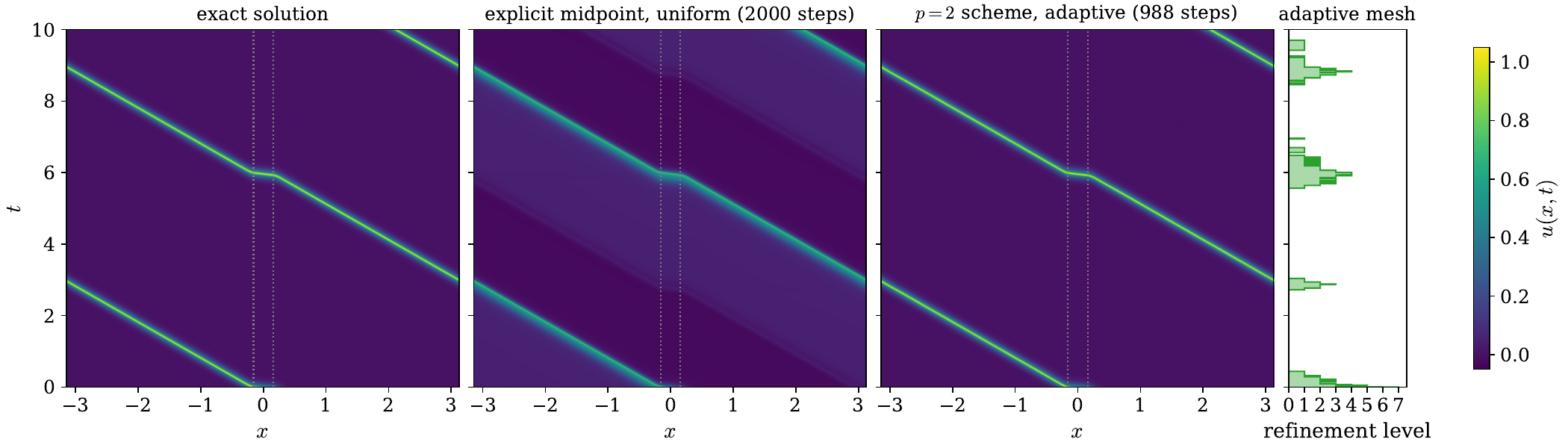}
\caption{\revision{Solution of the transport equation \eqref{eq:transport} with the localized coefficient spike shown in Figure~\ref{fig:err-coeff}, together with different numerical approximations.} The borders of the region where the coefficients $c$ spike are marked by dotted lines.}
\label{fig:spacetime-transport}
\end{figure}

\revision{The exact solution, the numerical approximation with $2000$ uniform explicit midpoint steps, and the adaptive approximation with $988$ time intervals are shown} in Figure~\ref{fig:spacetime-transport}. \revision{In the final plot of that figure}, we show the refinements chosen by the adaptive algorithm. The algorithm refines mostly when the initial mass interacts with the spike of the coefficients $c$. Some refinements are further chosen just in the region mirroring the region where the coefficients spike (\revision{i.e.,} are translated by half of the length of the interval $\pi$). This is not an artifact of the boundary, since the boundary is invisible to the parametrization (which enforces the periodic boundary conditions by setting the input nodes to $\sin(x+b_i)$ for fitted biases $b_1,\dots,b_{10}$). Instead, this seems to be a consequence of the parametric map $\Phi$ and the specific network architecture chosen for this experiment. 

\revision{As a result, the adaptive approximation is already free of visual artifacts and reaches an accuracy that requires roughly $10$ times as many steps with the explicit midpoint rule}, as is observed in the left plot of Figure~\ref{fig:err-coeff}. 

\subsection{\revision{Time-series latent ODE}}
\begin{figure}[h]
\begin{subfigure}{0.49\textwidth}
    \centering
    \includegraphics[width=\linewidth,height=0.679121\linewidth,trim=6bp 1bp 0bp 36bp,clip]{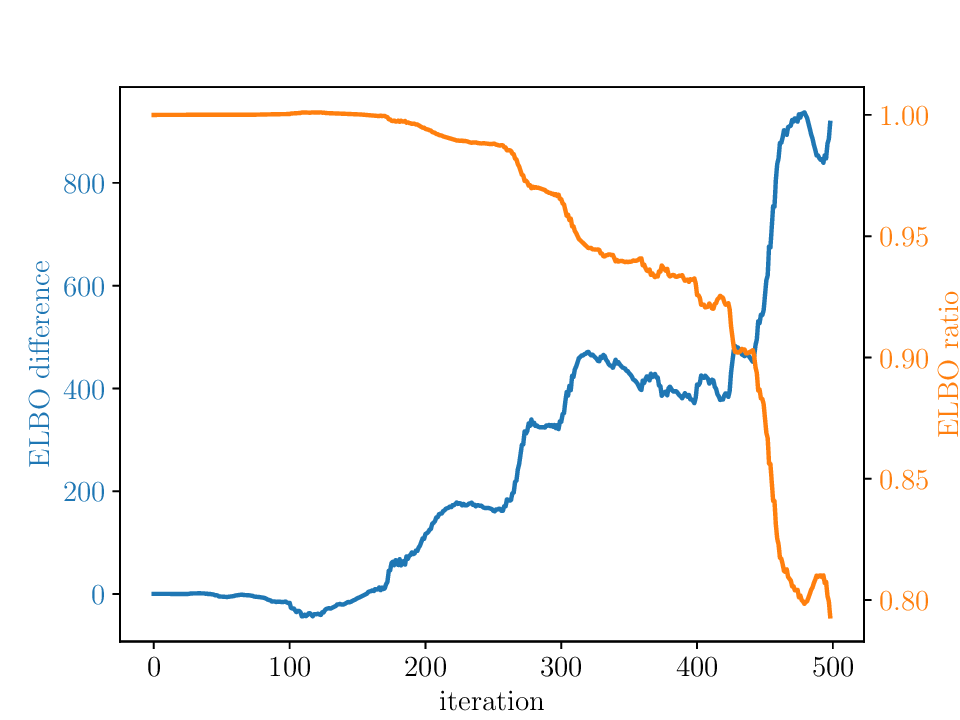}
    \caption{$n_{\rm max}=10$}
\end{subfigure}
\begin{subfigure}{0.49\textwidth}
    \centering
    \includegraphics[width=\linewidth,trim=6bp 1bp 0bp 36bp,clip]{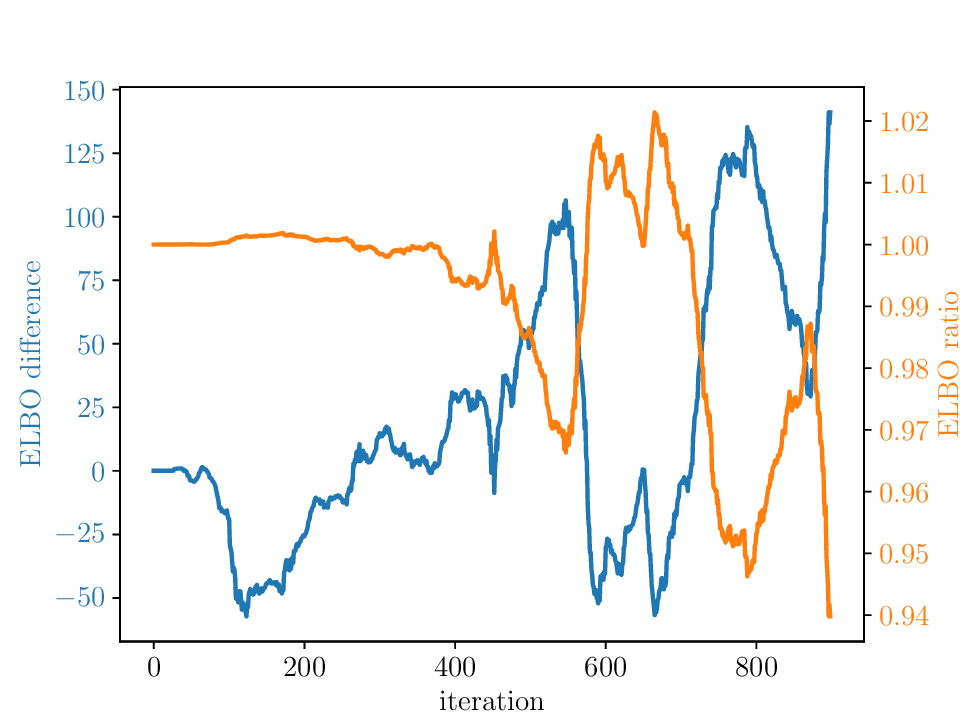}
    \caption{$n_{\rm max}=20$}
\end{subfigure}
\caption{\revision{Average ELBO loss versus the number of iterations. The solvers used} for inference and backpropagation use $n_{\rm max}=10$ in (a) and $n_{\rm max}=20$ in (b). \revision{Both the difference (adaptive--uniform) and the ratio (adaptive/uniform)} are depicted. Note that ELBO has a negative sign. $a=0.05$ in the definition of $\phi$.}
\label{fig:5nodes}
\end{figure}

\begin{figure}[h]
    \begin{subfigure}{0.49\textwidth}
    \centering
    \raisebox{-0.064675\linewidth}[0.788887\linewidth][0.064675\linewidth]{%
        \includegraphics[width=\linewidth,height=0.853562\linewidth,trim=13bp 1bp 41bp 38bp,clip]{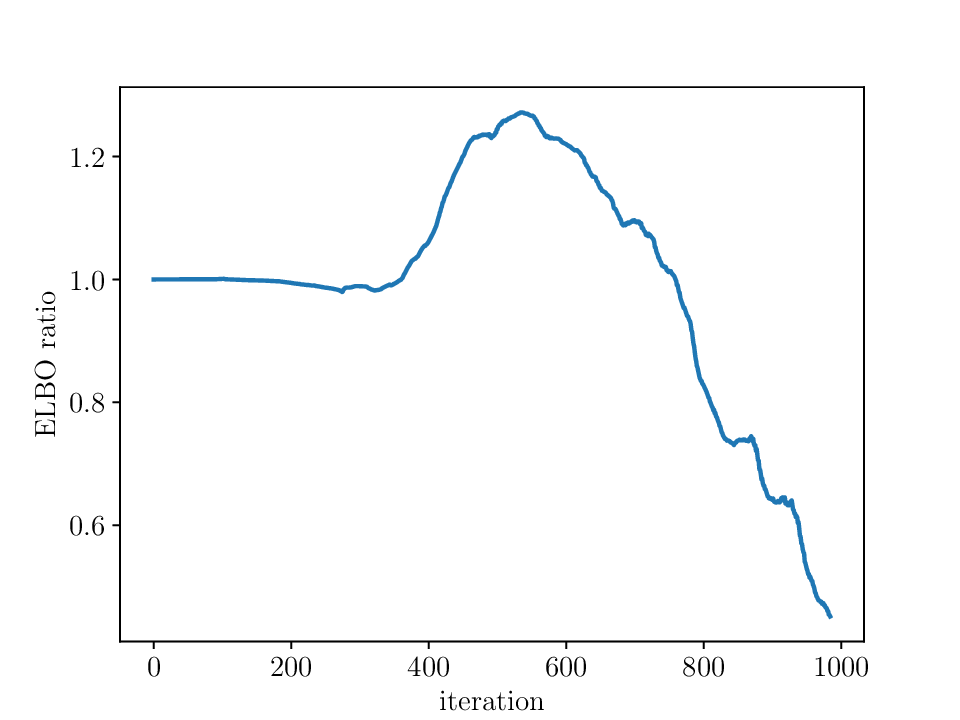}}
    \caption{adaptive/(non-adaptive)}
\end{subfigure}
\begin{subfigure}{0.49\textwidth}
    \centering
    \raisebox{0pt}[0.788887\linewidth][0.064675\linewidth]{%
        \includegraphics[width=\linewidth,trim=56bp 37bp 62bp 59bp,clip]{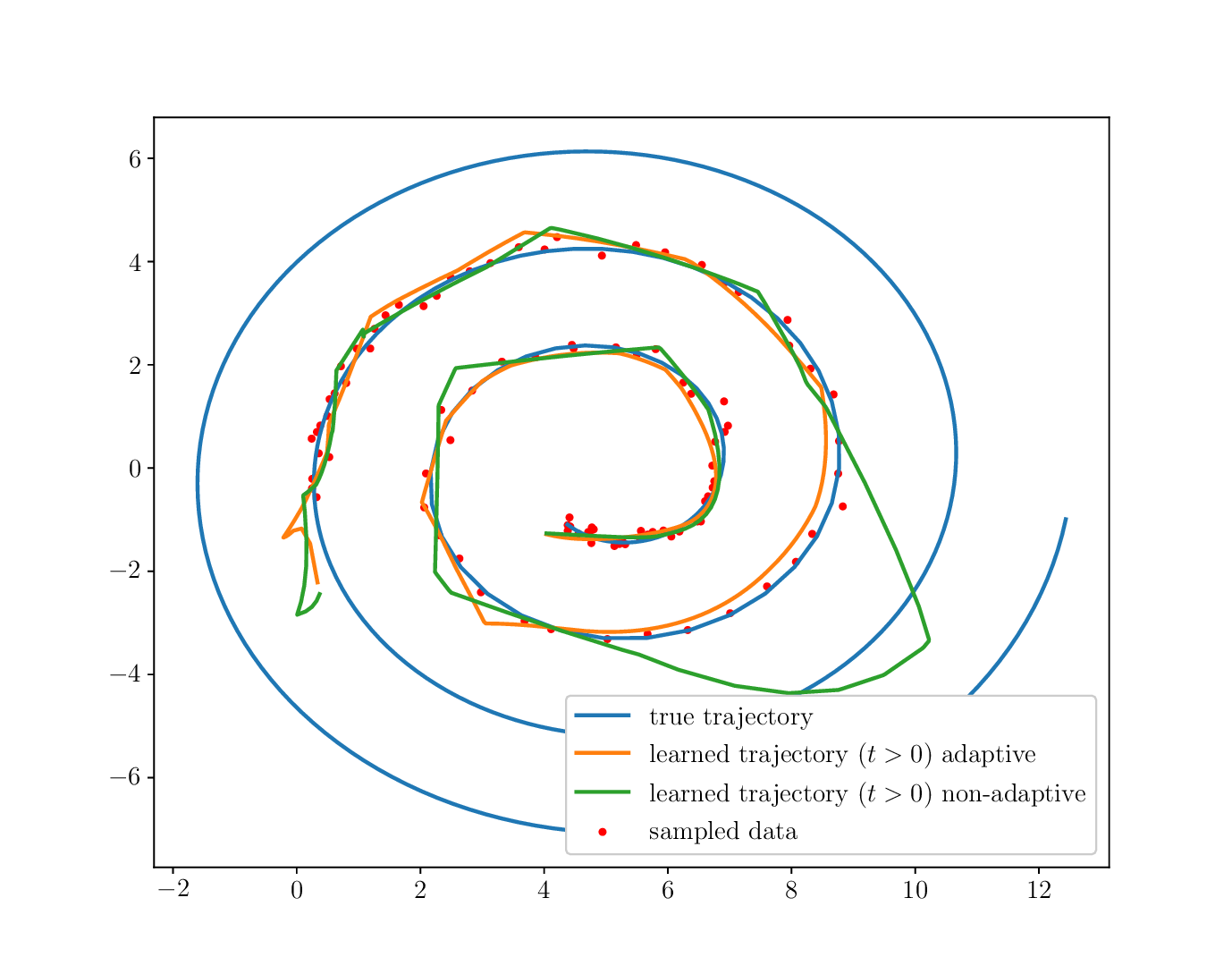}}
    \caption{Qualitative spiral reconstruction}
\end{subfigure}
\caption{(a): Performance of the adaptive training method applied to a spiral with $a=0.5$ in the definition of $\phi$.  (b): A qualitative picture comparing reconstructed spirals after training. \revision{Here, too, $a=0.5$.}}
\label{fig:qualitative}
\end{figure}

\begin{figure}[h]
    \begin{subfigure}{0.49\textwidth}
    \centering
    \includegraphics[width=\linewidth,trim=22bp 0bp 42bp 38bp,clip]{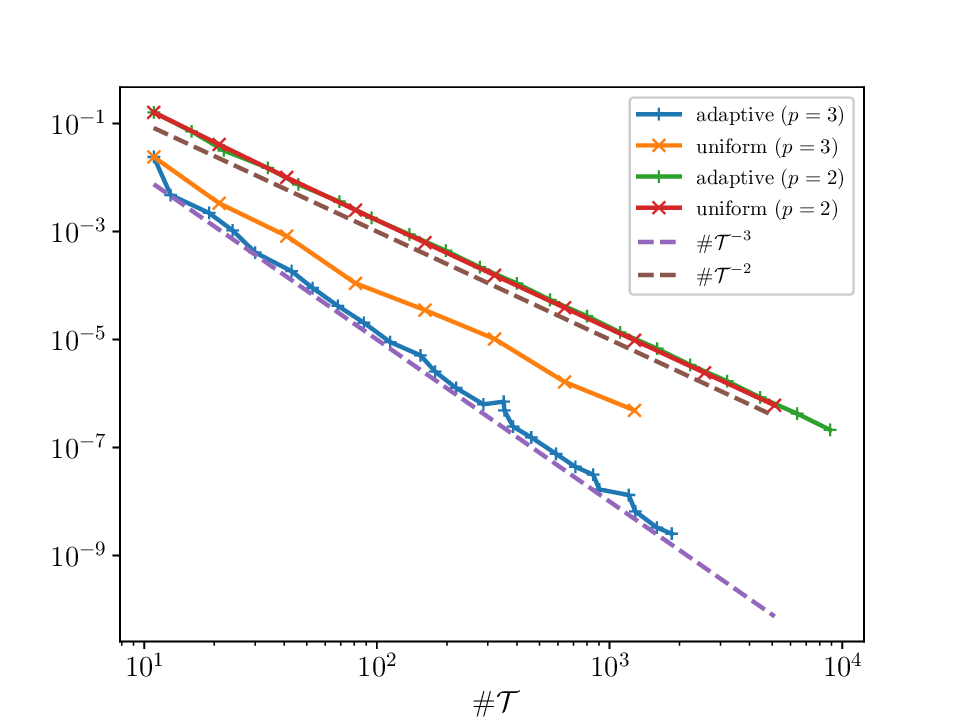}
    \caption{\revision{forward pass}}
\end{subfigure}
\begin{subfigure}{0.49\textwidth}
    \centering
    \includegraphics[width=\linewidth,trim=22bp 0bp 42bp 38bp,clip]{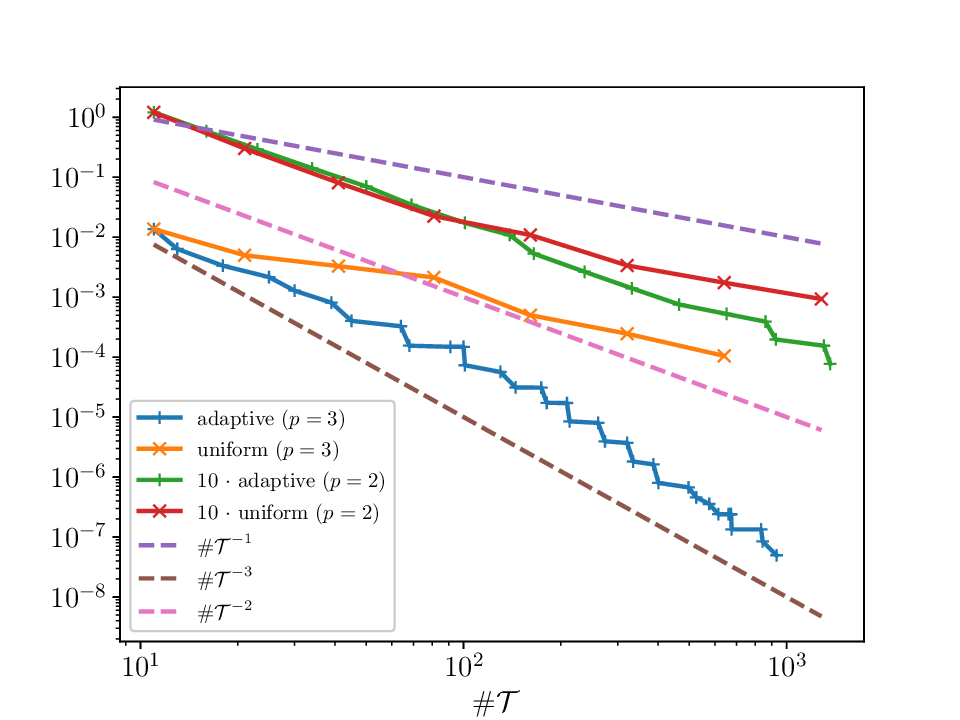}
    \caption{\revision{backward pass}}
\end{subfigure}
\caption{\revision{The explicit adaptive schemes of orders $p \in \{2,3\}$} with $\eta_\TT^\infty$ used for \revision{mesh refinement} \revision{applied to both the forward and backward passes} of the neural ODE. Note that in the backward pass, the estimator values for the second order scheme are multiplied by $10$ for better visibility. $\theta=0.5$ in all runs.}
 \label{fig:nn_p3}
\end{figure}

\revision{An example where a neural ODE network structure has several benefits is the reconstruction of dynamics from (irregular) samples in time, as laid out in} \cite[Section 5]{neuralODE}. The authors introduce a continuous-time generative approach to modeling time series and present an example which aims to reconstruct spirals, {\color{black} given in polar coordinates by
\begin{align*}
    r_1(t)=0.3\phi(t), \quad r_2(t)=0.3 \frac{50}{21-\phi(t)}, \quad
    \phi(t)=t, \quad t \in (0,20).
\end{align*}} We adapt this example by focusing on the spiral given by $r_1$, and adjusting it such that it is traversed very rapidly in a certain \revision{time window}. {\color{black} \revision{To achieve this, we choose}
\begin{align*}
    \phi(t)=t+5 \tanh\Big(\frac{t-10}{a}\Big), \quad a \in \{0.05, 0.5\}.
\end{align*}}This introduces a localized region of rapid variation in the observed trajectory around \(t=10\) and therefore \revision{provides a suitable example} for our adaptive procedure. We do not change the model architecture used in the paper. \revision{The loss function used to train the model is} the \textit{evidence lower bound (ELBO)}, commonly used in \textit{Variational Autoencoders} \cite{kingma2022autoencodingvariationalbayes}. \revision{First, we investigate} how adaptive algorithms can improve training of such latent ODE models. This is done by fixing a number of time steps $n_{\rm max}$ both for inference and backpropagation, \revision{and comparing the efficiency of training with adaptive and uniform methods} with $n_{\rm max}$ nodes. The results are shown in Figure~\ref{fig:5nodes} and Figure~\ref{fig:qualitative}. One observes advantages when using our adaptive solver compared to a \revision{fixed-grid method} of the same order. \revision{The advantage is more pronounced when a small value of $n_{\rm max}$ is prescribed.} In all experiments that investigate the training, the \revision{third-order} method from Figure~\ref{fig:rk} is used. We give a qualitative picture of the reconstruction abilities of both trained models in Figure~\ref{fig:qualitative}. The difference between adaptive and uniform refinement in training can also be observed there. After training, we investigated the convergence rate of the adaptive algorithm when applied to both the forward pass through the latent ODE and the backward pass, \revision{i.e.,} solving the adjoint problem \eqref{eq:ode_aug_simp}. We observe that for the forward pass, the adaptive \revision{third-order} method can retrieve its convergence order, while the uniform method \revision{only gives} a reduced order of convergence. This is even more pronounced in the backward pass, where both adaptive schemes ($p\in \{2,3\}$) converge with the optimal convergence rates, in stark contrast to uniform refinements. The results are shown in Figure~\ref{fig:nn_p3}. 
\section*{Declaration of AI use}
AI tools were used to assist in coding, identifying mistakes, typographical errors, inconsistencies and to improve the wording of the manuscript. The authors take full responsibility for its content.

\bibliographystyle{siamplain}
\bibliography{literature}

\begin{thebibliography}{10}

\bibitem{bdd}
{\sc P.~Binev, W.~Dahmen, and R.~DeVore}, {\em Adaptive finite element methods with convergence rates}, Numer. Math., 97 (2004), pp.~219--268.

\bibitem{axioms}
{\sc C.~Carstensen, M.~Feischl, M.~Page, and D.~Praetorius}, {\em Axioms of adaptivity}, Comput. Math. Appl., 67 (2014), pp.~1195--1253, \url{https://doi.org/10.1016/j.camwa.2013.12.003}, \url{http://dx.doi.org/10.1016/j.camwa.2013.12.003}.

\bibitem{ccgs13}
{\sc C.~Carstensen, D.~Gallistl, and M.~Schedensack}, {\em Adaptive nonconforming {C}rouzeix-{R}aviart {FEM} for eigenvalue problems}, Math. Comp., 84 (2015), pp.~1061--1087, \url{https://doi.org/10.1090/S0025-5718-2014-02894-9}, \url{https://doi.org/10.1090/S0025-5718-2014-02894-9}.

\bibitem{CCnonlin}
{\sc C.~Carstensen and N.~Nataraj}, {\em Adaptive {M}orley {FEM} for the von {K}\'{a}rm\'{a}n equations with optimal convergence rates}, SIAM J. Numer. Anal., 59 (2021), pp.~696--719, \url{https://doi.org/10.1137/20M1335613}, \url{https://doi.org/10.1137/20M1335613}.

\bibitem{ckns}
{\sc J.~M. Cascon, C.~Kreuzer, R.~H. Nochetto, and K.~G. Siebert}, {\em Quasi-optimal convergence rate for an adaptive finite element method}, SIAM J. Numer. Anal., 46 (2008), pp.~2524--2550.

\bibitem{neuralODE}
{\sc T.~Q. Chen, Y.~Rubanova, J.~Bettencourt, and D.~Duvenaud}, {\em Neural ordinary differential equations}, in Advances in Neural Information Processing Systems 31: Annual Conference on Neural Information Processing Systems 2018, NeurIPS 2018, December 3-8, 2018, Montr{\'{e}}al, Canada, S.~Bengio, H.~M. Wallach, H.~Larochelle, K.~Grauman, N.~Cesa{-}Bianchi, and R.~Garnett, eds., 2018, pp.~6572--6583, \url{https://proceedings.neurips.cc/paper/2018/hash/69386f6bb1dfed68692a24c8686939b9-Abstract.html}.

\bibitem{stevensonl2}
{\sc A.~Demlow and R.~Stevenson}, {\em Convergence and quasi-optimality of an adaptive finite element method for controlling {$L_2$} errors}, Numer. Math., 117 (2011), pp.~185--218, \url{https://doi.org/10.1007/s00211-010-0349-9}, \url{http://dx.doi.org/10.1007/s00211-010-0349-9}.

\bibitem{thetaopt}
{\sc L.~Diening and C.~Kreuzer}, {\em On the threshold condition for dörfler marking}, Computational Methods in Applied Mathematics, 21 (2021), pp.~557--567, \url{https://doi.org/doi:10.1515/cmam-2020-0041}, \url{https://doi.org/10.1515/cmam-2020-0041}.

\bibitem{instanceopt}
{\sc L.~Diening, C.~Kreuzer, and R.~Stevenson}, {\em Instance optimality of the adaptive maximum strategy}, Found. Comput. Math., 16 (2016), pp.~33--68, \url{https://doi.org/10.1007/s10208-014-9236-6}, \url{https://doi.org/10.1007/s10208-014-9236-6}.

\bibitem{diening2026optimalafemconvergencebulk}
{\sc L.~Diening and J.~Storn}, {\em Optimal afem convergence for any bulk parameter smaller than one}, 2026, \url{https://arxiv.org/abs/2610.39719}, \url{https://arxiv.org/abs/2610.39719}.

\bibitem{d1996}
{\sc W.~D{\"o}rfler}, {\em A convergent adaptive algorithm for {P}oisson's equation}, SIAM J. Numer. Anal., 33 (1996), pp.~1106--1124.

\bibitem{Ee}
{\sc W.~Ee}, {\em A proposal on machine learning via dynamical systems}, Communications in Mathematics and Statistics, 5 (2017), pp.~1--11, \url{https://doi.org/10.1007/s40304-017-0103-z}.

\bibitem{feischl2022}
{\sc M.~Feischl}, {\em Inf-sup stability implies quasi-orthogonality}, Math. Comp., 91 (2022), pp.~2059--2094, \url{https://doi.org/10.1090/mcom/3748}, \url{https://doi.org/10.1090/mcom/3748}.

\bibitem{nonsymm}
{\sc M.~Feischl, T.~F{\"u}hrer, and D.~Praetorius}, {\em {Adaptive {FEM} with optimal convergence rates for a certain class of nonsymmetric and possibly nonlinear problems}}, SIAM J. Numer. Anal., 52 (2014), pp.~601--625, \url{https://doi.org/10.1137/120897225}, \url{http://dx.doi.org/10.1137/120897225}.

\bibitem{regdyn}
{\sc M.~Feischl, C.~Lasser, C.~Lubich, and J.~Nick}, {\em Regularized dynamical parametric approximation}, Numer. Math., 158 (2026), pp.~2079--2123, \url{https://doi.org/10.1007/s00211-026-01551-5}, \url{https://doi.org/10.1007/s00211-026-01551-5}.

\bibitem{implicitadaptive}
{\sc M.~Feischl and D.~Niederkofler}, {\em Chapter nine - optimal adaptive implicit time stepping}, in Error Control, Adaptive Discretizations, and Applications, Part 4, F.~Chouly, S.~P. Bordas, R.~Becker, and P.~Omnes, eds., vol.~61 of Advances in Applied Mechanics, Elsevier, 2025, pp.~407--448, \url{https://doi.org/https://doi.org/10.1016/bs.aams.2025.09.001}, \url{https://www.sciencedirect.com/science/article/pii/S0065215625000250}.

\bibitem{Haber_2017}
{\sc E.~Haber and L.~Ruthotto}, {\em Stable architectures for deep neural networks}, Inverse Problems, 34 (2017), p.~014004, \url{https://doi.org/10.1088/1361-6420/aa9a90}, \url{http://dx.doi.org/10.1088/1361-6420/aa9a90}.

\bibitem{kingma2022autoencodingvariationalbayes}
{\sc D.~P. Kingma and M.~Welling}, {\em Auto-encoding variational bayes}, 2022, \url{https://arxiv.org/abs/1312.6114}, \url{https://arxiv.org/abs/1312.6114}.

\bibitem{li2026rateoptimalityafemmaximum}
{\sc Y.~Li}, {\em Rate optimality of afem with maximum marking}, 2026, \url{https://arxiv.org/abs/2611.03387}, \url{https://arxiv.org/abs/2611.03387}.

\bibitem{LN25}
{\sc C.~Lubich and J.~Nick}, {\em Regularized dynamical parametric approximation of stiff evolution problems}, arXiv preprint arXiv:2501.12118,  (2025).

\bibitem{explicit_a_posteriori}
{\sc H.~Ranocha and J.~Giesselmann}, {\em Stability of step size control based on a posteriori error estimates}, Computational Science and Engineering, 1 (2024), p.~1, \url{https://doi.org/10.1007/s44207-024-00001-0}, \url{https://doi.org/10.1007/s44207-024-00001-0}.

\bibitem{stevenson07}
{\sc R.~Stevenson}, {\em Optimality of a standard adaptive finite element method}, Found. Comput. Math., 7 (2007), pp.~245--269.

\bibitem{nonconf}
{\sc A.~Veeser and P.~Zanotti}, {\em Quasi-optimal nonconforming methods for symmetric elliptic problems. {I}---{A}bstract theory}, SIAM J. Numer. Anal., 56 (2018), pp.~1621--1642, \url{https://doi.org/10.1137/17M1116362}, \url{https://doi.org/10.1137/17M1116362}.

\bibitem{hairer_stiff}
{\sc G.~Wanner and E.~Hairer}, {\em Solving ordinary differential equations II}, vol.~375, Springer Berlin Heidelberg New York, 1996.

\end{thebibliography}
\end{document}